\documentclass[11pt]{article}
\usepackage[english]{babel}
\usepackage{amssymb, amsmath, amsthm} 
\usepackage[a4paper, centering]{geometry}
\usepackage{graphicx}
\usepackage{caption}
\usepackage{subfigure,color}
\usepackage[font={small,up}]{caption}

\numberwithin{equation}{section} 

\usepackage{pgfplots}
\usepackage{hyperref}
\usepackage{bm}

\usepackage{amsthm}
\makeatletter
\def\th@plain{%
  \thm@notefont{}% same as heading font
  \itshape % body font
}
\def\th@definition{%
  \thm@notefont{}% same as heading font
  \normalfont % body font
}
\makeatother
\newtheorem{thm}{Theorem}[section]
\newtheorem{prop}[thm]{Proposition}
\newtheorem{lem}[thm]{Lemma}
\newtheorem{cor}[thm]{Corollary}

\theoremstyle{definition}
\newtheorem{defn}[thm]{Definition}

\newtheorem{oss}[thm]{Remark}

\newcommand{\N}{\mathbb{N}}
\newcommand{\R}{\mathbb{R}}

\renewcommand{\epsilon}{\varepsilon}

\DeclareMathOperator{\dist}{dist}

\usepackage{enumitem}
\setlist[description]{nosep}

\title{Multiscale analysis of singularly perturbed finite dimensional gradient flows: the minimizing movement approach}

\author{{\scshape Giovanni Scilla}\\
Department of Mathematics and Applications ``R. Caccioppoli''\\ University of Naples ``Federico II''\\
Via Cintia, Monte S. Angelo - 80126 Naples \\
(ITALY)
 \\
\\
{\scshape Francesco Solombrino}\\
Department of Mathematics and Applications ``R. Caccioppoli''\\ University of Naples ``Federico II''\\
Via Cintia, Monte S. Angelo - 80126 Naples \\
(ITALY)}

\date{}

\begin{document}

\maketitle

\begin{abstract}
{We perform a convergence analysis of a discrete-in-time minimization scheme approximating a finite dimensional singularly perturbed gradient flow. We allow for different scalings between the viscosity parameter $\varepsilon$ and the time scale $\tau$.
 When the ratio $\frac{\epsilon}{\tau}$  diverges, we rigorously prove the convergence of this scheme to a (discontinuous) Balanced Viscosity solution of the quasistatic evolution problem obtained as formal limit, when $\varepsilon\to 0$,  of the gradient flow. We also characterize the limit evolution corresponding to an asymptotically finite ratio between the scales, which is of a different kind. In this case, a discrete interfacial energy is optimized at jump times}.
\end{abstract}

\noindent
{\bf Keywords:} gradient flow, singular perturbations, Balanced Viscosity solutions, variational methods, minimizing movement, rate-independent systems, crease energy\\
%\noindent
%{\bf Mathematics Subject Classification:}  %34K26, 34K18, 34C37, 47J35, 49J45

\section{Introduction}

This paper is concerned with the {convergence} analysis of the %limit behavior of 
minimizing movement approximations (also known as De Giorgi's approach to metric gradient flows, see, e.g.,~\cite{AmbGiSav}) of the singularly perturbed gradient flow system
\begin{equation}
\epsilon\dot{u}_\epsilon(t)=-\nabla_xF(t,u_\epsilon(t))\,,
\label{problem1}
\end{equation}
for  a time-depending driving energy $F$. The above system is a standard approximation for the quasistatic evolution problem
\begin{equation}\label{quasi-st}
\nabla_x F(t,u(t))=0.
\end{equation}
In its simplest formulation in Euclidean spaces (see also \cite[Section 3.2.1]{Braides} for some examples), the {minimizing movement approximation} method takes the form of a discrete-in-time variational scheme  
\begin{equation}
u^k_{\tau,\varepsilon}\in \text{argmin}\left\{F(t^k,x)+\frac{\varepsilon}{2\tau}\left\|x-u^{k-1}_{\tau,\varepsilon}\right\|^2,\quad x\in X\right\},\,\text{for } k=1,2,\dots,\lfloor T/\tau\rfloor,
\label{schemeintro}
\end{equation}
with initial point $u^0_{\tau,\varepsilon}=u^0$, where $\tau>0$ is a time scale, $\{t^k\}_k$ is the corresponding partition of $[0,T]$ defined by $t^k=k\tau$, $k=0,1,\dots,\lfloor T/\tau\rfloor$, and the term $\frac{\varepsilon}{2\tau}\left\|\,\cdot\,-u^{k-1}_{\tau,\varepsilon}\right\|^2$, containing the small \emph{viscosity parameter} $\epsilon$, penalizes the squared distance from the previous step $u^{k-1}_{\tau,\varepsilon}$. The role of this latter term in the above scheme is intuitively to avoid unnatural jumps, favoring {\it local} minimization of the energy instead of a global one. {Recursive schemes like \eqref{schemeintro} are well-known in Convex Optimization as \emph{proximal algorithms} (see, e.g., \cite[Ch. 5]{Bert}). Our focus will lie, however, on the case where $F$ is allowed to be \emph{nonconvex}, where discontinuous solutions to \eqref{quasi-st} are expected to arise in the limit as $\epsilon \to 0$  and $\tau \to 0$.}
\\
{\bf Our motivation.} {The scope of our analysis is twofold. On the one hand, we aim at complementing the abstract variational theory of gradient flows and of their numerical approximation. On the other hand, from a modeling perspective, in the limit as both $\epsilon \to 0$  and $\tau \to 0$ the scheme \eqref{schemeintro} provides a selection criterion for mechanically feasible solutions of \eqref{quasi-st}, whose properties are also investigated in our paper.}

{Concerning the first of the two aspects, we start by noticing that, under suitable smoothness of $F(t,x)$ and for small $\tau$, \eqref{schemeintro} is equivalent to an implicit Euler scheme (see, e.g., \cite{Bu}) for equation \eqref{problem1}, namely}
\begin{equation}
{\varepsilon\frac{u^k_{\tau,\varepsilon}-u^{k-1}_{\tau,\varepsilon}}{\tau}=-\nabla_xF(t^k, u^k_{\tau,\varepsilon}).}
\label{euler}
\end{equation}
{When $\varepsilon$ is fixed, under suitable assumptions on $F(t,x)$ the (interpolated) solutions $u_{\tau,\varepsilon}$ to \eqref{schemeintro} converge, as $\tau\to0$, to a function $u_\varepsilon(t)$ solving \eqref{problem1} for every $t\in[0,T]$ (see \cite[Section~213]{Bu} for details). In the limit $\varepsilon\to0$, one then expects to recover a solution $u(t)$ to \eqref{quasi-st}, which sits at local minima (on wells) of the potential energy $F(t,x)$, with the possible exception of a set of discontinuous points. The points where $u(t)$ is
discontinuous will be understood from \eqref{problem1} as the instantaneous and optimal transition from one stationary state to another at the onset of instability. A precise meaning to this heuristic picture has been given in the recent paper \cite{Ago-Rossi}, where the concept of Balanced Viscosity solutions to \eqref{quasi-st} has been introduced.  A first natural issue which we address in this paper is then taking the simultaneous limits $\epsilon \to 0$  and $\tau \to 0$ in equation \eqref{euler}. In such a case one expects that the sequential limit of taking $\tau \to 0$ first, followed by $\epsilon \to 0$, or the simultaneous limit $\epsilon ,\tau \to 0$ with $\tau/\epsilon \to 0$, will result in $u_{\tau,\varepsilon}$ converging to a Balanced Viscosity solution. This is shown to be true and is our first main result in Theorem \ref{main1}.}

{There is however another regime where one is interested in, and which we also analyse in this paper, that is when
$\tau \sim \lambda \varepsilon $ for a finite $\lambda \in (0, +\infty)$. A first reason for doing that comes from  numerical analysis. For the purposes of computational efficiency,  indeed, one typically does not want to have to take the time step $\tau$ to be extremely small. It is therefore a relevant issue to be addressed, whether one can choose a time step $\tau$ of order $\varepsilon$ and still obtain convergence to a meaningful evolution. Our second main result, Theorem \ref{main2}, positively answers this question. The limit evolution $u(t)$ satisfies again a local stability condition, which is in general stronger than \eqref{quasi-st}. Furthermore, an energy-dissipation balance is recovered, which accounts for the behavior at jumps in terms of a {\it discrete} optimization problem.}

{From a mechanics-oriented point of view, both Theorems  \ref{main1} and \ref{main2} provide an evolution $u(t)$ through {\it locally stable} states of the energy $F(t, \cdot)$, which, although being discontinuous, is {\it regulated}, i.e.\  left and right limits of $u(t)$ are well-defined at every time. The stored energy is decreased at discontinuity times, whereas the dissipation is optimal in an energetic sense. These are all relevant  features to be taken into account when dealing with discontinuous solutions to \eqref{quasi-st}, which are instead in general not satisfied by evolutions through  global minimizers of the energy (the so-called {\it energetic solutions}). These latter indeed show abrupt discontinuities, with no energy being dissipated along a jump, and are as a consequence very irregular.}

{While for the moment we will confine ourselves  to the case of smooth energy functionals in finite dimension, avoiding the additional technical issues coming from nonsmoothness and infinite dimensionality, extensions of our results in this direction can have relevant applications. We may mention for instance the experiments in \cite{Plos}, where transitions from the straight configuration to locally minimizing hemihelical configurations are observed under quasi-static releasing of a boundary load (we refer to \cite{Cicalese2017, Localmin} for a variational analysis of the related energy functional). We also remark that a forward Euler scheme, with a convex constraint in place of the penalization term and involving a time reparametrization, has been used in \cite{mnegri} to solve \eqref{quasi-st} and applied to a phase-field evolutionary model in brittle fracture. More in general}, let us mention that time-incremental minimization methods are by now a well-known tool in connection with the related setting of rate-independent evolutionary systems (see, e.g., \cite{MiThLev} and \cite{MiRou,Mielke} for an overview of the theory). There, viscous corrections have been often introduced (see, e.g., \cite{MiRoSav1,MiRoSav2,MiRoSav3, MinSav})  in order to better induce an evolution through local minimizers. In the setting of crack propagation, related approaches where the limit ratio $\frac{\varepsilon}{\tau}$ is kept finite have been considered (for instance in \cite{DalMasoToader,ACFS}).
The time-incremental minimization scheme for the viscous corrections to rate-independent systems carries however an additional dissipation term, taking (in the simplest case) the form  
\begin{equation}
u^k_{\tau,\varepsilon}\in \text{argmin}\left\{F(t^k,x)+\alpha\|x-u^{k-1}_{\tau,\varepsilon}\|+\frac{\varepsilon}{2\tau}\left\|x-u^{k-1}_{\tau,\varepsilon}\right\|^2,\quad x\in X\right\}\,,
\label{RIS}
\end{equation} 
with $\alpha>0$. Therefore, \eqref{schemeintro} can be seen as a degenerate case of \eqref{RIS}. {We however remark that the convergence analysis for \eqref{schemeintro} has to cope with relevant additional compactness issues in comparison with \eqref{RIS}. In this last case, one can profit of the term $\|x-u^{k-1}_{\tau,\varepsilon}\|$ which enforces $BV$-compactness for suitable interpolations of the scheme. Such an estimate is in general not available for \eqref{schemeintro} or even when dealing wih \eqref{problem1} in a time-continuous setting. On the other hand, the analysis of our problem is deeply inspired by recent results in the context of rate-independent systems, where, for instance, the notion of Balanced Viscosity solution has been originally introduced \cite{MiRoSav3}. Furthermore, the limit evolution provided by Theorem \ref{main2} (which, as we will discuss later, exhibits a sort of intermediate behavior between energetic and balanced viscosity solutions) can be considered a counterpart in our setting of the notion of {\it visco-energetic}  solutions to rate-independent systems introduced in  \cite{MinSav}.}

Before describing our results in detail, we review some recent contributions concerning \eqref{problem1} which serve as a starting point  for our analysis.
\\
{\bf Results present in literature.} A continuous theory about the limit evolutions of the solutions of the singularly perturbed problem (\ref{problem1}) in finite dimension has been developed by several authors (see, e.g.,~\cite{Zanini,Ago1,Ago-Rossi-Sav,Ago-Rossi,SciSol}). In particular, Agostiniani and Rossi in \cite{Ago-Rossi}, combining ideas from the variational approach to gradient flows with techniques for the vanishing viscosity approximation of the rate-independent systems, prove the existence of a limit curve $u$ that is a \emph{Balanced Viscosity solution} to \eqref{quasi-st}.
Along with  typical regularity, coercivity and power control assumptions on $F$ (which we also consider, see {\bf (F0)-(F2)}), a crucial role in the analysis is played by the assumption that the set of critical points $\mathcal{C}(t):=\left\{u\in X:\, \nabla_x F(t,u)=0\right\}$ consists of {\it isolated points} (also this one is assumed in  our paper, see {\bf(F3)} below). This allows for recovering the necessary  compactness through a careful analysis of the behavior at jumps. They indeed show that, up to a subsequence, the solutions $u_\epsilon$ of \eqref{problem1} pointwise converge, as $\epsilon\to0$, to a solution $u$ of the limit problem \eqref{quasi-st} defined at every $t\in[0,T]$ and such that:
\begin{enumerate}
\item[(1)] $u:[0,T]\to X$ is \emph{regulated}, i.e., the left and right limits $u_-(t)$ and $u_+(t)$ exist at every $t\in(0,T)$, and so do the limits $u_+(0)$ and $u_-(T)$;
\item[(2)] $u$ fulfills the \emph{energy balance}  
\begin{equation*}
F(t,u_+(t))+\mu([s,t])=F(s,u_-(s))+\int_s^t \partial_r F(r,u(r))\,dr\quad \forall\,0\leq s\leq t\leq T,
\end{equation*} 
where $\mu$ is a positive Radon measure with an at most countable jump set $J$;
\item[(3)] $u$ satisfies the stability condition
\begin{equation*}
\nabla_xF(t,u(t))=0
\end{equation*}
at all continuity points;
\item[(4)] $J$ coincides with the jump set of $u$ and the following jump relation holds:
\begin{equation*}
\mu(\{t\})=F(t,u_-(t))-F(t,u_+(t))=c_t(u_+(t),u_-(t)),\quad \forall\,t\in J,
\end{equation*}
\end{enumerate}
where the cost $c_t$ is defined as
\begin{equation}
c_t(u_1,u_2):=\inf\left\{\int_0^1\|\dot{\theta}(s)\|\|\nabla_x F(t, \theta(s))\|\,ds:\,\theta(0)=u_1, \theta(1)=u_2\right\}\,.
\label{costfunintro}
\end{equation}
In particular, transitions between (meta)stable states of the energy happen along (a finite union of) heteroclinic orbits of the unscaled gradient flow, an idea that has been already successfully used, in less generality, for instance in \cite{Zanini}.
Moreover, under additional assumptions (in the same spirit of our assumptions {\bf (F4)-(F5)}), the defect measure $\mu$ is  a pure jump measure, and its support coincides exactly with the jump set $J$ of $u$. {The auhors also provide some simple explicit examples \cite[Example 2.6]{Ago-Rossi} of balanced viscosity solutions which, as their heuristic meaning suggests, consist of a finite number of branches of local minimizers of $F(t, \cdot)$ parametrized by the time $t$, and jump exactly when reaching a degenerate  critical point. While energetic solutions jump as soon as possible, when global minimality gets lost,  balanced viscosity solutions jump instead as late as possible.}\\
{\bf Description of our results.} We now come to our results. We consider the discrete-in-time variational scheme  
\begin{equation*}
u^k_{\tau,\varepsilon}\in \text{argmin}\left\{F(t^k,x)+\frac{\varepsilon}{2\tau}\left\|x-u^{k-1}_{\tau,\varepsilon}\right\|^2,\quad x\in X\right\},\,\text{for } k=1,2,\dots,\lfloor T/\tau\rfloor,
\end{equation*}
with initial point $u^0_{\tau,\varepsilon}=u^0$, where $\tau>0$ is a time scale, $\{t^k\}_k$ is the corresponding partition of $[0,T]$ defined by $t^k=k\tau$, $k=0,1,\dots,\lfloor T/\tau\rfloor$, under the hypothesis that $X$ is a finite-dimensional Hilbert space, already considered in \cite{Ago-Rossi}. While the extension of our results to the infinite-dimensional setting carries additional technical issues, which we do not deal with in the present paper, we believe that the techniques we develop here already provide the necessary insight into the problem and are likely to be significant tools also in infinite dimension.

{A first mathematical issues to be overcome is recovering  a compactness theorem for suitable interpolations of the sequence $\{u^k_{\tau,\varepsilon}\}$. This is a relevant point, since in general we expect convergence to a discontinuous solution, whose total variation, furthermore, can not be a priori estimated. A different strategy, involving some nontrivial refinement of the arguments in \cite{Ago-Rossi}, has then to be pursued.  For this, we first introduce the the piecewise constant $\tilde{u}_{\epsilon,\tau}$ and the piecewise affine $\bar{u}_{\epsilon,\tau}$ interpolations of $\{u^k_{\tau,\varepsilon}\}$, respectively. Some elementary energy bounds are then obtained, which however make use of both these interpolations (see \eqref{enestimateimproved}). Hence, we can not directly apply the compactness arguments of \cite[Theorem 1]{Ago-Rossi}.  A first step is using the energy bounds to show that, if the ratio  $\epsilon/\tau$ stays bounded, tthe mismatch between $\tilde{u}_{\epsilon,\tau}$ and $\bar{u}_{\epsilon,\tau}$ is vanishing in $L^2$, and, whenever $\tau <\!\!< \epsilon$, even in $L^\infty$ (Corollary \ref{samelimit-}).  Since we are interested in taking a pointwise limit at all times $t$ and we also allow for $\tau \sim \epsilon$, this is still not enough to recover an everywhere defined limit function. A crucial point , which is typical of the discrete setting, is then to show that}, when doing a discrete transition between metastable states, the gradients $\nabla_x F(t^k, u^k)$ stay bounded away from zero. In terms of the interpolations, this amounts to requiring that, when  the gradient of the piecewise affine interpolations is bounded away from zero, then also the gradient of the piecewise constant ones is.
This fact is established in Lemma \ref{boundGrad} as a consequence of  {\bf(F3)} and of the Euler-Lagrange conditions \eqref{parallel}. It will allow us to show that transitions  always happen at the price of a strictly positive cost, and therefore only a countable number of them is allowed. We also prove that, if $u$ is a pointwise limit of $\tilde{u}_{\epsilon,\tau}$, the limit stored energy $s\to F(s,u(s))$ is a function of bounded variation on $[0,T]$ (Proposition~\ref{equi-boundvariation}). After establishing compactness, this is useful to recover existence of the one-sided limits at discontinuity points.

With Theorem~\ref{compactnessthm} we show the aforementioned compactness property for $(\tilde{u}_{\epsilon,\tau}(t))$: under the assumption that the ratio $\epsilon/\tau$ stays bounded, we first pass to the limit (along a subsequence independent of $t$) on $(\tilde{u}_{\epsilon,\tau}(t))_{\epsilon,\tau}$, and prove that they converge for all $t$ to a regulated function $u(t)$. This limit function in general satisfies the stability condition
\begin{equation}
\nabla_x F(t,u(t))=0
\label{stazintro2}
\end{equation}
for all $t \in [0,T]\setminus J$, where $J$ is the (at most countable) jump set of $u$. When the ratio $\epsilon/\tau$ tends to a finite limit, a stronger form  of the stability conditions (Proposition~\ref{improve}), more suitable for this regime, can be also deduced. Furthermore, we  use \eqref{stazintro2} to eventually show that indeed also the piecewise affine interpolations $(\bar{u}_{\epsilon,\tau}(t))_{\epsilon,\tau}$ converge to the same limit $u(t)$ at all its continuity points (Corollary \ref{samelimit}). 

After that compactness and stability properties of the limit evolution have been established, we show that the limit evolution $u(t)$ satisfies, independently of the limit ratio $\frac{\epsilon}{\tau}$, a balance between the stored energy and  the power spent  along the evolution in an interval of time $[s, t]\subset [0, T]$, up to a positive dissipation cost which is concentrated on the jump set of $u$, or equivalently on the jump set of the energy $t\to F(t, u(t))$. Namely, we prove in Theorem~\ref{genbalance} that there exists a {\it positive atomic measure} $\mu$, with ${\rm supp}(\mu)=J$, such that
\begin{equation*}
F(t,u_+(t))+\mu([s,t])=F(s,u_-(s))+\int_s^t \partial_r F(r,u(r))\,dr,
\end{equation*} 
for all $0\leq s\leq t\leq T$. In order to do this, we borrow some ideas from \cite[Section 7]{MinSav}. While the ``$\,\le\,$'' inequality above can be obtained by passing to the limit in the considered minimization scheme, the ``$\,\ge\,$'' ensues from the stability condition, and requires the additional assumptions {\bf(F4)-(F5)} on the energy, which are instead not necessary in order to recover compactness (see Section \ref{generals}).

In order to characterize the dissipation cost, we have to leave the unified framework valid until now and to distinguish between the regimes $\tau<\!<\!\epsilon$ and $\tau\sim\epsilon$. In the first case (Section \ref{casoinfinito}), we retrieve the notion of {\it Balanced Viscosity solution} introduced in the continuum setting in \cite{Ago-Rossi}, as we show (Proposition~\ref{costotrovato}) that
\begin{equation*}
\mu(\{t\})=c_t(u_+(t),u_-(t)),\quad \forall\,t\in J\,.
\end{equation*}
where $c_t$ coincides with the cost (\ref{costfunintro}). {Also here, the proof strategy  has  to face some difficulties that are peculiar of the discrete setting. In particular, to be in a position to apply the lower semicontinuity properties of the cost $c_t$ discussed in \cite{Ago-Rossi}, we have to show that the integral terms
\[
\int\|\nabla_x F(s, \tilde{u}_{\epsilon,\tau}(s))\|\|\dot{\bar{u}}_{\epsilon,\tau}(s)\|\,ds 
\]
appearing in the energy estimates, and
\[
\int \|\nabla_x F(s, \bar{u}_{\epsilon,\tau}(s))\|\|\dot{\bar{u}}_{\epsilon,\tau}(s)\|\,ds\,,
\]
which instead only contains the piecewise affine interpolations $\bar{u}_{\epsilon,\tau}$, carry an asymptotically equivalent dissipation at discontinuity times. While we already know that, under the assumption $\tau<\!\!<\epsilon$,  the mismatch between the two interpolations is uniformly small by \eqref{uniforme}, this still might be not enough to establish the desidered equivalence, since no $L^1$-estimate for $\|\dot{\bar{u}}_{\epsilon,\tau}(s)\|$ is at hand. A recursive construction exploiting {\bf(F3)} is then needed to show that an optimal decrease of the energy can be realised at discontinuities via a {\it finite} number of transitions between metastable states, and the integral
\[
\int \|\dot{\bar{u}}_{\epsilon,\tau}(s)\|\,ds
\]
can only blow up near the endpoints of these transitions, where $\nabla_x F(\cdot, \cdot)$ is small. Combining with \eqref{uniforme}, this eventually allows us to show in Theorem \ref{main1} that, if  $\tau<\!\!<\epsilon$, the limit $u(t)$ is exactly a Balanced Viscosity solution.}

In the regime $\tau\sim\epsilon$ (Section \ref{casofinito}), {for $\lambda$ being the limit ratio of $\frac{\epsilon}{\tau}$} we introduce a parametric cost $c^\lambda$ as {an interfacial energy, solving a discrete optimization problem (see, e.g.,~\cite{BC,SciVal})}. Namely, for $\mathcal{R}_\lambda$ being the \emph{residual stability function} defined by \eqref{residual}, we set
\begin{equation}
c^\lambda(t,u,v):=\displaystyle\inf\left\{\sum_{i=0}^{N-1}\frac{\lambda}{2}\|w^i-w^{i+1}\|^2+\sum_{i=0}^{N}\mathcal{R}_\lambda(t,w^i)\right\},
\label{costintro}
\end{equation}
where  the infimum is taken over all families $(w^i)_{i\in \mathbb{N}}$ with $w^0=u$ and $w^N=v$ and all $N \in \mathbb{N}$. The function $\mathcal{R}_\lambda$ provides a measure of the failure of the \emph{stability condition} for $(t,u)\in[0,T]\times X$. Roughly speaking, the jump transitions between $u_-(t)$ and $u_+(t)$ are described by discrete trajectories $W=(w^i)$ such that each $w^{i+1}$ is a minimizer of the incremental time scheme with fixed $t$ and datum $w^i$, where the first point of the chain is $w^0=u_-(t)$ and the last one is $w^N=u_+(t)$. A more general form of \eqref{costintro}, with sums parametrized by a continuous parameter, has been considered in~\cite[Definition~3.5]{MinSav} for the rate-independent setting. The necessity of taking sums on (possibly non discrete) compact subsets of $\mathbb{R}$ was motivated there with the fact that an infinite number of transitions could in general occur at discontinuity times. In our case, we are able to exclude this with assumption {\bf(F3)} and an inductive construction, and we can work  with the simpler definition \eqref{costintro}. {In particular, a crucial lower semicontinuity property for the cost $c^\lambda$ is recovered in Theorem \ref{semicont}}. In our second main result Theorem \ref{main2} we completely characterize the limit evolution in the regime $\tau\sim\epsilon$, showing the equality
\begin{equation}\label{cost++}
\mu(\{t\})=c^\lambda(t,u_-(t),u_+(t))\,.
\end{equation}
{As we already mentioned, the solution concept provided by Theorem \ref{main2} is in some sense intermediate between energetic and balanced viscosity solutions. Indeed, $u(t)$ has to satisfy the stability property}
\begin{equation}\label{sloc}
F(t,u(t))\le F(t,v)+\frac\lambda2\|v-u(t)\|^2
\end{equation}
{for all $v\in X$, which is in general stronger than \eqref{quasi-st}, but still retains a local character since the square distance from $u(t)$ is penalised. As an effect, solutions jump later than in the energetic case, with a strictly positive dissipation obeying \eqref{cost++}, but in general before reaching a degenerate critical point as in the balanced viscosity case. For a simple one-dimensional example with $X=\mathbb{R}$ one  may think of the double well energy $F(t,u)=\frac14(u^2-1)^2-(t-1)u$ for $t\in [0,2]$. The energetic solution follows the lower branch of the curve $u^3-u=t-1$ for $t\le 1$, when it jumps from $u_-=-1$ to $u_+=1$ with no energy being dissipated. The balanced viscosity solutions, instead, continues following the lower branch till it reaches, at $t=1+\frac2{3\sqrt3}$, the degenerate critical point $u=-\frac1{\sqrt3}$ and then jumps to the upper branch. If one takes, instead, for instance $\lambda=\frac12$, it follows from \eqref{sloc} that $t=1+\frac2{3\sqrt3}$ is a continuity point for the evolution provided by  Theorem \ref{main2}, with $u=\frac{2}{\sqrt3}$ being the only possible value. Therefore $u(t)$ has  jumped before, well inside the time interval $(1, 1+\frac2{3\sqrt3})$.} 
\\
{\bf Plan of the paper.} The paper is organized as follows. In Section~\ref{preliminar} we fix notation and provide some preliminaries, recalling the main assumptions on the energies $F(t,x)$ we will adopt throughout {and also commenting on their generic character.} Section~\ref{minimizingmovements} deals with the time incremental minimization scheme for the gradient flow system, whose solution is a discrete-in-time sequence $\{u^k\}$. As customary in this setting, we introduce suitable interpolations of this values (piecewise constant/piecewise affine) and prove some basic inequalities.
In Section~\ref{stability}, we prove  compactness of the interpolations and stability properties of the limit evolution.
In Section~\ref{energybalance} we show that the limit evolution $u$ fulfills an energy balance with a cost concentrated on the jump points (Theorem~\ref{genbalance}), independently of the limit of the ratio $\epsilon/\tau$, at the price of additional, but still general, assumptions {\bf (F4)-(F5)}.
Only afterwards, to characterize the dissipation cost, we will  distinguish between the  two cases when $\epsilon/\tau\to+\infty$ (Section~\ref{casoinfinito}), and $\epsilon/\tau$ is bounded (Section~\ref{casofinito}). The two main results of the paper, containing complete characterizations of the limit evolutions depending on the limit ratio  $\epsilon/\tau$ are stated in Theorems \ref{main1}, and \ref{main2}, respectively.

\section{Preliminaries and assumption on the energy}\label{preliminar}
Preliminarily, let us fix some general notation that will be used throughout. We aim at describing quasistatic evolutions driven by  a time-dependent, possibly nonconvex energy functional $F: [0, T]\times X\to\R$, with $T>0$.
Throughout the paper we assume that $(X,\|\cdot\|)$ is a Euclidean space with dimension $n\geq1$, endowed with inner product $\langle \cdot, \cdot \rangle$.  Given $x \in X$ and $M > 0$, we will denote by $B(x,M)$ the closed ball centered at $x$ with radius $M$. When the ball is centered at $0$, the shortcut $B_M$ will be used.  

We now recall the definition of {\it regulated} functions, which will play an important role in the sequel.

\begin{defn}
A function $u:[0,T]\to X$ is said to be \emph{regulated} if for each $s\in[0,T]$ there exist the one-sided limits $u_+(s)$ and $u_-(s)$ defined as
\begin{equation*}
u_+(s):=\lim_{h\to0^+} u(s+h),
\end{equation*}
and
\begin{equation*}
u_-(s):=\lim_{h\to0^-} u(s+h).
\end{equation*}
\end{defn}

The existence of the above limits immediately implies that, for each $N \in \mathbb{N}$, the set of points $t$ where $\|u_+(t)-u_-(t)\|\ge \frac1N$ cannot have accumulation points. It follows that the jump set of  a regulated function is at most countable. In particular, $u_+$ is a right-continuous Lebesgue representative of $u$ and $u_-$ is a left-continuous one. 

It is well-known that a function of bounded variation  $f\in {\rm BV}([a,b];\R)$ is a real-valued regulated function. The representatives $f_+$ and $f_-$ are in this case good representatives in the sense of \cite[Theorem 3.28]{AFP}: as shown there, the distributional derivative $Df$ (which is a Radon measure) satisfies 
\begin{equation}\label{teoremafond}
Df([s,t])=f_+(t)-f_-(s)
\end{equation}
for any $s,t\in[a,b]$ with $s\le t$.

\subsection{Assumptions on the energy}
We require that the energy functional satisfy the same assumptions considered in \cite{Ago-Rossi}, namely, 

\begin{enumerate}
\item[({\bf F0)}] $F\in C^1([0,T]\times X)$,

\item[{\bf (F1)}] denoting with $\mathcal{F}$ the map $\mathcal{F}(u):=\sup_{t\in[0, T]}|F(t,u)|$, it holds 
for all $\rho>0$ that the sublevel set $\{u\in X:\, \mathcal{F}(u)\leq \rho\}$ is bounded;

\item[{\bf(F2)}] there exist $C_1,C_2>0$ such that 
\begin{equation*}
 |\partial_t F(t,u)|\leq C_1 (F(t,u)+C_2),\quad \text{for every }\, (t,u)\in[0,T]\times X,
\end{equation*}
where $\partial_t F$ denotes the partial derivative of $F(t,x)$ with respect to $t$;
\item[{\bf (F3)}] for any $t\in[0,T]$, the set of critical points
\begin{equation}
\mathcal{C}(t):=\Bigl\{u\in X:\, \nabla_x F(t,u)=0\Bigr\}
\label{criticalset}
\end{equation}
where $\nabla_x F$ denotes the gradient of $F(t,x)$ with respect to $x$, consists of isolated points. In particular, the set $\mathcal{C}(t)\cap B$ is finite whenever $B\subset X$ is relatively compact.
\end{enumerate}

We note that from {\bf (F2)} and the Gronwall's inequality we get
{\begin{equation}
F(t,u)\leq \left(F(s,u)+C_2\right)e^{C_1(t-s)}-C_2,
\label{gronwall}
\end{equation}}
for all $0\le s\le t\le T$ and $u\in X$.

The above assumptions will be enough to establish compactness and stability of the limit functions in Section \ref{stability}, while in Section \ref{energybalance} we will be forced to consider two additional assumptions {\bf(F4)}-{\bf(F5)} in order to recover an energy balance.
We will namely assume that the driving energy $F(t,x)$ satisfies:

\begin{enumerate} 
\item[{\bf (F4)}] For any $t\in[0,T]$ and for any $u\in \mathcal{C}(t)$, 
\begin{equation*}
\displaystyle \mathop{\lim\inf}_{v\to u}\frac{F(t,v)-F(t,u)}{\|\nabla_xF(t,v)\|}\geq0.
\end{equation*}
\item[{\bf (F5)}] For  fixed $u \in X$, the function $t\to \nabla_xF(t,u)$ is Lipschitz continuous on $[0,T]$, locally uniformly w.r.t. $u$. 
\end{enumerate}
Condition {\bf(F5)} is satisfied, e.g., if the mixed derivative $\partial_t\nabla_xF(\cdot,\cdot)$ is continuous on $[0,T]\times X$.

{Before concluding this section, we briefly discuss the generic character of assumptions {\bf (F3)-(F4)}, which are indeed satisfied by a very broad class of potentials. Since nondegenerate stationary points of smooth functionals are always isolated by the Implicit Function theorem, condition {\bf (F3)} has usually only to be checked at points $(t,u)$ where $u$ is a \emph{degenerate} critical point of $F(t, \cdot)$; i.e., $\nabla_x^2F(t,u)$ is non-invertible. It has been shown in \cite{Ago-Rossi-Sav} that energies $F\in C^3([0,T]\times X)$, satisfying the so-called \emph{transversality conditions} (see, e.g., \cite[Def.~6.1]{Ago-Rossi} for a precise definition), comply with {\bf (F3)}. These conditions are well-known in the realm of bifurcation theory: essentially, they imply that a curve $\phi(t)$ of critical points for $F(t, \cdot)$ has a fold whenever crossing a degenerate one. Remarkably, \cite[Thm.6.3]{Ago-Rossi-Sav} proves that, for every given energy functional $F(t,u)$, {\it any arbitrarily small quadratic perturbation} thereof (up to a meagre set in a suitable topology) fulfills the transversality conditions, which are therefore generic in a rigorous mathematical sense.  
The paper \cite{Ago-Rossi-Sav} also provides extensions to infinite-dimensional examples, such as a prototypical integral energy functional in elasticity, say
\begin{equation}\label{elas}
F(t,u):=\int_\Omega \left(\frac12|\nabla u(x)|^2+\mathcal{W}(u(x))-\ell(t,x)u(x)\right)\, dx,
\end{equation}
where $\Omega\subset\mathbb{R}^n$ is a connected reference configuration, $X:=H^2(\Omega)\cap H^1_0(\Omega)$ is  the space of admissible deformations $u$, $\ell\in C^4([0,T];L^2(\Omega))$ is a smooth applied load, and $\mathcal{W}(u):=(u^2-1)^2/4$ is a double-well potential. In \cite[Example 3.8]{Ago-Rossi-Sav} it is shown that the above functional, whose critical points solve a semilinear elliptic equation, complies with the transversality conditions up to possibly adding an arbitrarily small quadratic perturbation, chosen out of a meagre set in a suitable topology.}

{We also remark that even in the case of bifurcating branches  of critical points from a trivial state, where transversality conditions are not satisfied, assumption {\bf (F3)} still holds, for instance, at a bifurcation from simple eigenvalues (see \cite{CrRa} and, e.g., \cite[Proposition 4.3]{Localmin} for a recent application to a Kirchhoff rod model). }

{Condition {\bf(F4)} easily holds, for instance, if $F(t, \cdot)$ is convex for fixed $t$, since in this case $\nabla_xF(t,u)=0$ implies $F(t,v)\geq F(t,u)$ for every $v\in X$.  In the case of nonconvex potentials, {\bf (F4)} can be directly deduced (see \cite[Remark 2.5]{Ago-Rossi}) whenever $F(t,\cdot)$ complies with the celebrated \L ojasiewicz inequality; namely, for every $(t,u)\in[0,T]\times\mathcal{C}(t)$, there exist $\theta\in(0,1)$ and $C,R>0$ such that for every $v\in B_R(u)$ it holds
\begin{equation}\label{loja}
|F(t,v)-F(t,u)|^\theta\leq C\|\nabla_x F(t,u)\|\,.
\end{equation}
If $X$ has finite dimension, the above inequality is always satisfied whenever, at each fixed $t$, $F(t, \cdot)$ is an analytic function of $u$ (see \cite{Lo}). Hence, finite elements discretizations of \eqref{elas} fulfill {\bf (F4)}.  More in general, the  \L ojasiewicz inequality is  strictly related to the geometric notion of \emph{subanalytic} function (see \cite[Def.~6.4]{Ago-Rossi}). For the model case
\begin{equation}\label{model}
F(t,u)=E(u)-\langle \ell(t), u \rangle
\end{equation} 
where $E$ is a stored energy and $\ell(t)$ an applied load, it follows from \cite[Thm.~3.1]{BAL} that \eqref{loja} (and therefore {\bf (F4)}) is satisfied at every $(t,u)\in[0,T]\times\mathcal{C}(t)$ whenever $E\colon X \to \mathbb{R}$ is subanalytic. We also notice, that for $F$ as in \eqref{model}, condition {\bf (F5)} reduces to require that $\ell \colon [0,T]\to X$ is Lipschitz continuous.
For a relevant application of the \L ojasiewicz inequality to nonlinear evolution equations and gradient flows in infinite dimensions, we refer to \cite{Simon}. }

\section{The minimizing movement approach}\label{minimizingmovements}

Let $\tau>0$ be a time step. We consider a partition $\Pi_\tau:=\{t^k\}_k$ of the time-interval $[0,T]$ defined by $t^k:=k\tau$, where $k=0,1,\dots,m:=\lfloor T/\tau\rfloor$ and $\lfloor x\rfloor$ denotes the integer part of $x$. We note that $|\Pi_\tau|=m+1$ and we set $\mathcal{K}_\tau:=\{1,\dots,m\}$. We construct a recursive sequence $\{u^k_{\tau,\varepsilon}\}_k$ by solving the iterated minimum problem
\begin{equation}
u^k_{\tau,\varepsilon}\in \text{argmin}\left\{F(t^k,x)+\frac{\varepsilon}{2\tau}\left\|x-u^{k-1}_{\tau,\varepsilon}\right\|^2,\quad x\in X\right\},\,\text{for } k\in\mathcal{K}_\tau,
\label{scheme}
\end{equation}
with initial point $u^0_{\tau,\varepsilon}=u^0$.  
In order to simplify the notation, we will drop the dependence on $\tau,\varepsilon$ in $u^k_{\tau,\varepsilon}$ and we will denote it by $u^k$.

\begin{prop}\label{intmin}
Let {\bf (F0)-(F2)} hold. Then the minimum problem {\rm(\ref{scheme})} has solutions, and for $k\in\mathcal{K}_\tau$ we have the following upper energy estimate
\begin{equation}
F(t^k,u^k)-F(t^{k-1},u^{k-1})+\frac{\varepsilon}{2\tau}\|u^k-u^{k-1}\|^2\leq\int_{t^{k-1}}^{t^k}\partial_r F(r,u^{k-1})\,dr.
\label{stimaa}
\end{equation}
\end{prop}

\proof
With given $u^{k-1}$, any minimizer $x$ in (\ref{scheme}) at step $k$ satisfies the minimality condition
\begin{equation}
F(t^k,x)+\frac{\varepsilon}{2\tau}\|x-u^{k-1}\|^2\leq F(t^k,u^{k-1}),
\label{compare}
\end{equation}
since $u^{k-1}$ is a competitor in the minimization procedure. Setting $\mathcal{F}_k(x):=F(t^k,x)$, (\ref{compare}) implies that in order to determine $u^k$ we should minimize the lower semi-continuous functional $\mathcal{F}_k$ on the compact sublevel $F(t^k,\cdot)\leq F(t^k,u^{k-1})$. Thus, such an $u^k$ exists. Now we show (\ref{stimaa}). Rewriting (\ref{compare}) for $x=u^k$ and subtracting $F(t^{k-1},u^{k-1})$ from both the sides we immediately get
\begin{equation*}
\begin{split}
F(t^k,u^{k})-F(t^{k-1},u^{k-1})+\frac{\varepsilon}{2\tau}\|u^k-u^{k-1}\|^2&\leq F(t^k,u^{k-1})-F(t^{k-1},u^{k-1})\\
                                                                                               &= \int_{t^{k-1}}^{t^k}\partial_r F(r,u^{k-1})\,dr.
\end{split}
\end{equation*}
\endproof

\begin{prop}[\emph{A priori} estimate]
It holds
\begin{equation}
F(t^k,u^k)\leq(F(0,u^0)+C_2)e^{C_1 t^k}-C_2, \quad \text{for every }k\in\mathcal{K}_\tau.
\label{gronwa}
\end{equation}
\end{prop}

\proof
From {\bf (F2)}, (\ref{gronwall}) and (\ref{stimaa}) we obtain
\begin{equation}
\begin{split}
F(t^k,u^k)&\leq F(t^{k-1},u^{k-1}) + (F(t^{k-1},u^{k-1})+C_2)(e^{C_1(t^k-t^{k-1})}-1)\\
             &=(F(t^{k-1},u^{k-1})+C_2)e^{C_1(t^k-t^{k-1})}-C_2,
\end{split}
\label{stima3.7}
\end{equation}
and iterating over $k$ the left hand side of (\ref{stima3.7}) we finally get
\begin{equation*}
F(t^k,u^k)+C_2\leq(F(0,u^0)+C_2)\displaystyle\prod_{l=1}^ke^{C_1(t^l-t^{l-1})}=(F(0,u^0)+C_2)e^{C_1 t^k}, \quad \text{for }k\in\mathcal{K}_\tau.
\end{equation*}
\endproof

We denote by $\tilde{u}_{\epsilon,\tau}$ the left-continuous piecewise constant interpolant of $\{u^k\}_k$, defined by 
\begin{equation}
\tilde{u}_{\epsilon,\tau}(t):=u^{k},\,\text{ for } t\in(t^{k-1},t^k],\, k\in\mathcal{K}_\tau, 
\label{piecewiseconstant}
\end{equation}
and $\tilde{u}_{\epsilon,\tau}(T)=u^{m}$, and by $\bar{u}_{\epsilon,\tau}$ the piecewise affine interpolation of $\{u^k\}_k$, defined by
\begin{equation}
\bar{u}_{\epsilon,\tau}(t):=\frac{u^k-u^{k-1}}{\tau}(t-t^{k-1})+u^{k-1},\quad t\in(t^{k-1},t^k].
\label{piecewiseaffine}
\end{equation}
We note that
\begin{equation}
\bar{u}_{\epsilon,\tau}(t^k)=\tilde{u}_{\epsilon,\tau}(t^k),\, \text{for every }\, t^k\in\Pi_\tau.
\label{identity}
\end{equation}

For any $s\in[0,T]$, we denote by $\hat{s}_\tau$ the least point of partition $\Pi_\tau$ which is greater or equal to $s$; i.e., it is defined as
\begin{equation}
\hat{s}_\tau:=\min\{r\in\Pi_\tau:\, r\geq s\}.
\label{node}
\end{equation}
Note that $\hat{s}_\tau\searrow s$ as $\tau\to0$.

We also define a piecewise constant interpolation of the values $F(t^k,\cdot)$, setting
\begin{equation*}
F_\tau(t,u):= F(t^k,u),\quad \text{if }t\in(t^{k-1},t^k],\, \text{ for every }\, u\in X.
\end{equation*}

\begin{oss}

{From assumption {\bf(F0)} we deduce that, in the limit as $\tau \to 0$,
\begin{equation}
F_\tau(t,u)\to F(t,u)\quad \text{ and }\quad\nabla_x F_\tau(t,u)\to \nabla_xF(t,u),
\label{(*)}
\end{equation}
uniformly with respect to $(t,u) \in [0,T]\times B_M$.} 
\end{oss}
We note that by {\bf (F2)} and (\ref{gronwa}) there exists a positive uniform constant $C$ such that
\begin{equation*}
\mathcal{F}(\tilde{u}_{\epsilon,\tau})=\sup_{t\in[0,T]}|F(t, \tilde{u}_{\epsilon,\tau}(t))|\leq C,
\end{equation*}
and then, by assumption {\bf (F1)} on the compactness of $\{\mathcal{F}(\cdot)\leq C\}$, there exists $M>0$, independent from $\epsilon$ and $\tau$, such that
\begin{equation}
\|\tilde{u}_{\epsilon,\tau}(t)\|\leq M,
\label{compactfunc}
\end{equation}
for all $t\in[0,T]$. We summarize all the previous observations in the following Corollary. The energy estimate \eqref{energyestimate} that we derive here is not sharp, as we will prove in Section~\ref{energybalance}. However, it will be enough to obtain the main compactness result stated by Theorem~\ref{compactnessthm}.

\begin{cor}
Under assumptions {\bf (F0)-(F2)}, the functions $\tilde{u}_{\epsilon,\tau}$ and $\bar{u}_{\epsilon,\tau}$ satisfy the following estimates:
\begin{enumerate}
\item[\rm(i)] for every $s,t\in[0,T]$, $s<t$, it holds the energy estimate
\begin{equation}
\begin{split}
F_\tau(t, \tilde{u}_{\epsilon,\tau}(t))&+\frac{\varepsilon}{2}\int_{\hat{s}_\tau}^{\hat{t}_\tau}\|\dot{\bar{u}}_{\epsilon,\tau}(r)\|^2\,dr\\
&\leq F_\tau(s, \tilde{u}_{\epsilon,\tau}(s))+\int_{\hat{s}_\tau}^{\hat{t}_\tau}\partial_r F(r,\tilde{u}_{\epsilon,\tau}(r))\,dr.
\end{split}
\label{energyestimate}
\end{equation}
\item[\rm(ii)] for every $t\in[0,T]$, it holds the a priori estimate
\begin{equation*}
F_\tau(t,\tilde{u}_{\epsilon,\tau}(t))\leq C_0e^{C_1 \hat{t}_\tau}-C_2,
\end{equation*}
where $C_0:=F(0,u^0)+C_2$.
\end{enumerate}
\end{cor}

\begin{oss}
For any $s,t\in[0,T]$, $s<t$,  \eqref{energyestimate} can be also rewritten in the equivalent form
\begin{equation}
\begin{split}
F_\tau(t, \tilde{u}_{\epsilon,\tau}(t))&+\frac{1}{4}\int_{\hat{s}_\tau}^{\hat{t}_\tau}\left(\frac{1}{\varepsilon}\|\nabla_xF_\tau(r,\tilde{u}_{\epsilon,\tau}(r))\|^2+{\epsilon}\|\dot{\bar{u}}_{\epsilon,\tau}(r)\|^2\right)\,dr\\
&\leq F_\tau(s, \tilde{u}_{\epsilon,\tau}(s))+\int_{\hat{s}_\tau}^{\hat{t}_\tau}\partial_r F(r,\tilde{u}_{\epsilon,\tau}(r))\,dr.
\end{split}
\label{enestimateimproved}
\end{equation}

Indeed, if we write the Euler-Lagrange equations for the minimum problem \eqref{scheme}, we get
\begin{equation}
\frac{\varepsilon}{\tau}(u^k-u^{k-1})=-\nabla_xF(t^k,u^k),
\label{parallel}
\end{equation}
from which testing by $u^k-u^{k-1}$ and using the Cauchy inequality we obtain
\begin{equation*}
\begin{split}
\frac{\varepsilon}{2}\frac{\|u^k-u^{k-1}\|^2}{\tau}&=\frac{1}{2}\langle-\nabla_xF(t^k,u^k),{u^k-u^{k-1}}\rangle\\
&=\frac{\tau}{4\varepsilon}\|\nabla_xF(t^k,u^k)\|^2+\frac{\epsilon}{4\tau}\|u^k-u^{k-1}\|^2,
\end{split}
\end{equation*}
whence 
\begin{equation*}
\begin{split}
&\frac{\varepsilon}{2}\int_{\hat{s}_\tau}^{\hat{t}_\tau}\|\dot{\bar{u}}_{\epsilon,\tau}(r)\|^2\,dr\\
&=\frac{1}{4}\int_{\hat{s}_\tau}^{\hat{t}_\tau}\left(\frac{1}{\varepsilon}\|\nabla_xF_\tau(r,\tilde{u}_{\epsilon,\tau}(r))\|^2+{\epsilon}\|\dot{\bar{u}}_{\epsilon,\tau}(r)\|^2\right)\,dr,
\end{split}
\end{equation*}
which, together with \eqref{energyestimate}, immediately implies \eqref{enestimateimproved}.
\end{oss}

Endowed with the previous result, we can prove that the mismatch between the  interpolations $\tilde{u}_{\epsilon,\tau}$ and $\bar{u}_{\epsilon,\tau}$ is vanishing in $L^2$, and, whenever $\tau <\!\!< \epsilon$, even in $L^\infty$. We start with a result in the discrete setting.

\begin{prop} Let $\{u^k\}_k$ be defined by {\rm(\ref{scheme})}. Then
\begin{enumerate}
\item[{\rm(i)}] if $\epsilon,\tau\to0$ and there exists $C>0$ such that $\frac{\tau}{\epsilon}\leq C$, then
\begin{equation*}
\sum_{k\in\mathcal{K}_\tau}\tau\|u^{k}-u^{k-1}\|^2\to0;
\end{equation*}
\item[{\rm(ii)}] if $\epsilon,\tau\to0$ are such that $\frac{\tau}{\epsilon}\to0$ then
\begin{equation*}
\|u^{k}-u^{k-1}\|\to0
\end{equation*}
uniformly with respect to $k\in\mathcal{K}_\tau$.
\end{enumerate}
\label{tend0}
\end{prop}
\proof
(i) We have
\begin{equation*}
\begin{split}
\sum_{k\in\mathcal{K}_\tau}\tau\|u^{k}-u^{k-1}\|^2&=\sum_{k\in\mathcal{K}_\tau}\tau^2\int_{t^{k-1}}^{t^k}\|\dot{\bar{u}}_{\epsilon,\tau}(r)\|^2\,dr=\sum_{k\in\mathcal{K}_\tau}\frac{\tau^2}{\epsilon}\epsilon\int_{t^{k-1}}^{t^k}\|\dot{\bar{u}}_{\epsilon,\tau}(r)\|^2\,dr\\
&=\tau\frac{\tau}{\epsilon}\left(\epsilon\int_0^T\|\dot{\bar{u}}_{\epsilon,\tau}(r)\|^2\,dr\right)\leq\widetilde{C}\tau\to0,
\end{split}
\end{equation*}
where we used \eqref{energyestimate}.\\
(ii) We have
\begin{equation*}
\begin{split}
\|u^{k}-u^{k-1}\|&=\frac{\sqrt{\tau}}{\sqrt{\epsilon}}\frac{\sqrt{\epsilon}}{\sqrt{\tau}}\|u^{k}-u^{k-1}\|\leq\sqrt{\frac{{\tau}}{{\epsilon}}}\left(\epsilon\int_{t^{k-1}}^{t^k}\|\dot{\bar{u}}_{\epsilon,\tau}(r)\|^2\,dr\right)^\frac{1}{2}\\
&\leq \widetilde{C}^\frac{1}{2}\sqrt{\frac{{\tau}}{{\epsilon}}}\to0,
\end{split}
\end{equation*}
where we used \eqref{energyestimate} again.

\endproof

The previous proposition may be re-read in terms of the interpolations as stated below.

\begin{cor}
Let $\tilde{u}_{\epsilon,\tau}$ and $\bar{u}_{\epsilon,\tau}$ be defined as in \eqref{piecewiseconstant} and \eqref{piecewiseaffine}, respectively. Then
\begin{enumerate}
\item[\rm(i)] if $\epsilon,\tau\to0$ and there exists $C>0$ such that $\frac{\tau}{\epsilon}\leq C$, then
\begin{equation}\label{l2conv}
\int_0^T\|\tilde{u}_{\epsilon,\tau}(r)-\bar{u}_{\epsilon,\tau}(r)\|^2\,dr\to0;
\end{equation}
\item[\rm(ii)] if $\epsilon,\tau\to0$ are such that $\frac{\tau}{\epsilon}\to0$, then 
\begin{equation}
\tilde{u}_{\epsilon,\tau}(t)-\bar{u}_{\epsilon,\tau}(t)\to0
\end{equation}
uniformly with respect to $t\in[0,T]$.
\end{enumerate}
\label{samelimit-}
\end{cor}

\proof
Let $t\in(t^{k-1},t^k]$ for some $t^{k-1},t^k\in\Pi_\tau$. Then by a direct computation we have the estimate
\begin{equation}\label{confrontz}
\|\tilde{u}_{\epsilon,\tau}(t)-\bar{u}_{\epsilon,\tau}(t)\|=\frac{\|u^k-u^{k-1}\|}{\tau}|\tau-t+t^{k-1}|\leq\|u^k-u^{k-1}\|,
\end{equation}
from which (ii) follows by Proposition~\ref{tend0}(ii). Concerning (i), we note that
\begin{equation*}
\int_0^T\|\tilde{u}_{\epsilon,\tau}(r)-\bar{u}_{\epsilon,\tau}(r)\|^2\,dr=\sum_{k\in\mathcal{K}_\tau}\int_{t^{k-1}}^{t^k}\frac{\|u^{k}-u^{k-1}\|^2}{\tau^2}(\tau-t+t^{k-1})^2\leq\sum_{k\in\mathcal{K}_\tau}\tau\|u^{k}-u^{k-1}\|^2,
\end{equation*}
and the thesis follows by Proposition~\ref{tend0}(i).
\endproof

Since by {\bf (F0)} $F\in C^1([0,T]\times X)$, corresponding to the compact convex subset $B_M$ there exists a concave modulus of continuity $\omega: [0,+\infty)\to[0,+\infty)$ such that $\displaystyle\lim_{\beta\searrow0}\omega(\beta)=\omega(0)=0$ and
\begin{equation*}
\Bigl|\partial_r F(t,u)-\partial_r F(t,v)\Bigr|\leq\omega(\|u-v\|),\quad \text{ for every }\, u,v\in B_M,
\end{equation*}
uniformly with respect to $t\in[0,T]$.

The following proposition, very useful in the sequel, provides a connection between the integral term involving $\partial_r F(r, u^{k-1})$ contained in the upper energy estimate (\ref{stimaa}) and the corresponding one obtained when we replace $u^{k-1}$ by $\tilde{u}_{\epsilon,\tau}(r)=u^k$, showing that the global error committed when summing over the points of the partition $\Pi_\tau$ is vanishing as $\epsilon\to0$.

\begin{prop} It holds
\begin{equation*}
\sum_{k\in\mathcal{K}_\tau} \left|\int_{t^{k-1}}^{t^k}\Bigl[\partial_r F(r, u^{k-1})-\partial_r F(r, u^{k})\Bigr]\,dr\right|\to0
\end{equation*}
as $\epsilon\to0$.
\label{convergence1}
\end{prop}

\proof
Let $\omega$ be a (concave) modulus of continuity for $\partial_r F(t, u)$ on $[0,T]\times B_M$. We have
\begin{equation*}
\begin{split}
\sum_{k\in\mathcal{K}_\tau} \left|\int_{t^{k-1}}^{t^k}\Bigl[\partial_r F(r, u^{k-1})-\partial_r F(r, u^{k})\Bigr]\,dr\right|&\leq T\sum_{k\in\mathcal{K}_\tau}\frac{\tau}{T}\omega(\|u^k-u^{k-1}\|)\\
& \leq T\omega\left(\frac{1}{T}\sum_{k\in\mathcal{K}_\tau}\tau\|u^k-u^{k-1}\|\right),
\end{split}
\end{equation*}
where the latter tends to 0 since by Cauchy inequality we get
\begin{equation*}
\sum_{k\in\mathcal{K}_\tau} \tau \|u^k-u^{k-1}\|\leq \left(\sum_{k\in\mathcal{K}_\tau} \tau\right)^{\frac{1}{2}}\left(\sum_{k\in\mathcal{K}_\tau} \tau\|u^k-u^{k-1}\|^2\right)^{\frac{1}{2}},
\end{equation*}
and the first term in the right hand side is bounded, while the latter term tends to 0 by Proposition~\ref{tend0}(i).
\endproof

An analogous result holds when replacing the piecewise constant interpolations $\tilde{u}_{\epsilon,\tau}$ by the piecewise affine interpolations $\bar{u}_{\epsilon,\tau}$, as proved in the following proposition.

\begin{prop}
\begin{equation*}
\int_0^T \Bigl|\partial_r F(r, \bar{u}_{\epsilon,\tau}(r))-\partial_r F(r, \tilde{u}_{\epsilon,\tau}(r))\Bigr|\,dr\to0
\end{equation*}
as $\epsilon\to0$.
\end{prop}
\proof
We have the estimate
\begin{equation*}
\begin{split}
&\int_0^T \Bigl|\partial_r F(r, \bar{u}_{\epsilon,\tau}(r))-\partial_r F(r, \tilde{u}_{\epsilon,\tau}(r))\Bigr|\,dr\\
&=\sum_{k\in\mathcal{K}_\tau} \int_{t^{k-1}}^{t^k}\Bigl|\partial_r F(r, \bar{u}_{\epsilon,\tau}(r))-\partial_r F(r, \tilde{u}_{\epsilon,\tau}(r))\Bigr|\,dr\\
&\leq \sum_{k\in\mathcal{K}_\tau} \int_{t^{k-1}}^{t^k}\omega\left(\frac{\|u^k-u^{k-1}\|}{\tau}|t-t^{k-1}-\tau|\right)\\
&\leq \sum_{k\in\mathcal{K}_\tau} \tau \omega(\|u^k-u^{k-1}\|)\leq T\omega\left(\frac{1}{T}\sum_{k\in\mathcal{K}_\tau}\tau\|u^k-u^{k-1}\|\right),
\end{split}
\end{equation*}
where we used also the monotonicity properties of $\omega$. The conclusion easily follows with an analogous argument as in Proposition~\ref{convergence1}.
\endproof

Another technical tool useful in the proof of the main result of compactness (Theorem~\ref{compactnessthm}) is the following lemma, that provides a comparison between the gradient of $F$ computed at the piecewise affine interpolations $\bar{u}_{\epsilon,\tau}$ and the gradient of $F_\tau$ computed along the piecewise constant interpolations $\tilde{u}_{\epsilon,\tau}$. More precisely, also $\|\nabla_xF_\tau(t,\tilde{u}_{\epsilon,\tau}(t))\|$ is bounded away from zero on the set where $\|\nabla_xF(t,\bar{u}_{\epsilon,\tau}(t))\|$ is.

\begin{lem}\label{boundGrad}
Let $\epsilon,\tau\to0$ be sequences for which there exists $C>0$ such that $\frac{\tau}{\varepsilon}\leq C$, and consider the interpolants $\tilde{u}_{\epsilon,\tau}$ and $\bar{u}_{\epsilon,\tau}$ be defined as in \eqref{piecewiseconstant} and \eqref{piecewiseaffine}, respectively. Let $\omega$ be a uniform modulus of continuity for $\nabla_x F(\cdot,\cdot)$ on $[0,T]\times B_M$, with $M$ given by \eqref{compactfunc}. Then, setting
\begin{equation*}
\alpha:=\min_{t\in[a,b]}\|\nabla_xF(t,\bar{u}_{\epsilon,\tau}(t))\|
\end{equation*}
we have, for any $\beta>0$ such that $\omega(C\beta)+\omega(\tau)+\beta<\alpha$,
\begin{equation*}
\|\nabla_xF_\tau(t,\tilde{u}_{\epsilon,\tau}(t))\|\geq\beta,\,\text{ for every }\, t\in[a,b].
\end{equation*}
\label{lemmafogli}
\end{lem}
\proof
We argue by contradiction and assume that there exists some $t\in[a,b]$ such that
\begin{equation}
\|\nabla_xF_\tau(t,\tilde{u}_{\epsilon,\tau}(t))\|<\beta.
\label{bou}
\end{equation}
Let $t^k\in\Pi_\tau$ be such that $\nabla_xF_\tau(t,\tilde{u}_{\epsilon,\tau}(t))=\nabla_xF(t^k,u^k)$. Then from \eqref{parallel} and \eqref{bou} we get
\begin{equation*}
\|u^k-u^{k-1}\|=\frac{\tau}{\epsilon}\|\nabla_xF(t^k,u^k)\|<C\beta.
\end{equation*}
Correspondingly to the same $t$, we also obtain
\begin{equation*}
\begin{split}
\|\nabla_xF_\tau(t,\bar{u}_{\epsilon,\tau}(t))\| &= \|\nabla_xF(t^k,\bar{u}_{\epsilon,\tau}(t))\|\\
&\leq \|\nabla_xF(t^k,\bar{u}_{\epsilon,\tau}(t))-\nabla_xF(t^k,u^k)\|+\|\nabla_xF(t^k,u^k)\|\\
&\leq \omega(\|u^k-u^{k-1}\|) +\beta\leq \omega(C\beta) +\beta.
\end{split}
\end{equation*}
Now,
\begin{equation*}
\begin{split}
\|\nabla_xF(t,\bar{u}_{\epsilon,\tau}(t))\|&\leq \|\nabla_xF(t,\bar{u}_{\epsilon,\tau}(t))-\nabla_xF_\tau(t,\bar{u}_{\epsilon,\tau}(t))\|+\|\nabla_xF_\tau(t,\bar{u}_{\epsilon,\tau}(t))\|\\
&\leq \omega(\tau)+\omega(C\beta) +\beta,
\end{split}
\end{equation*}
so that taking the infimum of the left hand side with respect to $t\in[a,b]$ we get a contradiction.
\endproof

\begin{prop}
Let $\tilde{u}_{\epsilon, \tau}:[0,T]\to X$ be defined as in {\rm(\ref{piecewiseconstant})}. Assume that there exists a function $u:[0,T]\to X$ such that $\tilde{u}_\epsilon(s)\to u(s)$ for every $s\in[0,T]$ as $(\epsilon, \tau)\to(0, 0)$. Then the function 
\begin{equation*}
s\to f(s):= F(s,u(s))-\int_0^s \partial_r F(r,u(r))\,dr
\end{equation*}
is nonincreasing. In particular, $s\to F(s,u(s))$ is of bounded variation on $[0,T]$.
\label{equi-boundvariation}
\end{prop}

\proof
For every $s\in[0,T]$, we denote by $\hat{s}_\tau$ the corresponding node of the partition $\Pi_\tau$ defined as in (\ref{node}). We set
\begin{equation*}
f_{\epsilon, \tau}(s):= F_\tau(s,\tilde{u}_{\epsilon, \tau}(s))-\int_0^{\hat{s}_\tau}\partial_r F(r,\tilde{u}_{\epsilon, \tau}(r))\,dr.
\end{equation*} 
Let $0\leq s_1\leq s_2\leq T$. If we introduce the notation $k_{i\tau}:=\hat{s}_{i\tau}/\tau$, $i=1,2$, we then have
\begin{equation*}
\begin{split}
f_{\epsilon, \tau}(s_1)-f_{\epsilon, \tau}(s_2)&= F_\tau(s_1,\tilde{u}_{\epsilon, \tau}(s_1))-F_\tau(s_2,\tilde{u}_{\epsilon, \tau}(s_2))+\int_{\hat{s}_{1\tau}}^{\hat{s}_{2\tau}}\partial_r F(r,\tilde{u}_{\epsilon, \tau}(r))\,dr\\
& =o(1)+ \sum_{k=k_{1\tau}+1}^{k_{2\tau}}\left(F(t^{k-1},u^{k-1})-F(t^k,u^k)+\int_{t^{k-1}}^{t^k}\partial_r F(r,u^{k-1})\,dr\right)\\
&\ge o(1)
\end{split}
\end{equation*}
as $(\epsilon,\tau)\to(0,0)$, where we used Proposition~\ref{convergence1} and since  each summand is nonnegative by \eqref{stimaa}. In the limit we get $f(s_1)\geq f(s_2)$. The rest of the statement follows by the absolute continuity of the integral.
\endproof

\section{Compactness and stability properties}\label{stability}

Our first main result of this paper is the following compactness result, which we prove under the only assumption that the ratio $\epsilon/\tau$ stays bounded.
Within this general assumption, according to \eqref{l2conv}  the limits of the interpolations  $(\tilde{u}_{\epsilon,\tau}(t))_{\epsilon,\tau}$ and  $(\bar{u}_{\epsilon,\tau}(t))_{\epsilon,\tau}$ coincide in general only up to a null set. Since instead we are interested in {\it pointwise} convergence, we will then first pass to the limit (along a subsequence independent of $t$) on $(\tilde{u}_{\epsilon,\tau}(t))_{\epsilon,\tau}$, and prove that they converge for all $t$ to a regulated function $u(t)$. This limit function satisfies the stability condition
\begin{equation}
\nabla_x F(t,u(t))=0
\label{staz}
\end{equation}
for all $t \in [0,T]\setminus J$, where this latter is the (at most countable) jump set of $u$. We will then use \eqref{staz} to show that indeed \eqref{l2conv} can be refined, and also $(\bar{u}_{\epsilon,\tau}(t))_{\epsilon,\tau}$ converge to $u(t)$ at  all continuity points. In particular, the left- and right- continuous representatives $u_+(t)$ and $u_-(t)$, in terms of which both the stability condition \eqref{stazionario} and the energy balances \eqref{ebalance} and \eqref{ebalance2} are formulated, are determined independently of the interpolations.
We finally also show that \eqref{staz} can be improved to a stronger stability condition whenever $\tau\sim \epsilon$ (see Proposition \ref{improve} below). 

We first state and prove our compactness result.

\begin{thm} Assume that {\bf (F0)-(F3)} hold and let $\tilde{u}_{\epsilon,\tau}: [0,T]\to X$ be the piecewise constant functions defined in \eqref{piecewiseconstant}, interpolating the discrete solutions of the minimum problem \eqref{scheme}, with $\tilde{u}_{\epsilon,\tau}(0)=u^0$. Then all the sequences $(\varepsilon_j,\tau_j)_{j\in\mathbb{N}}$ satisfying $(\varepsilon_j,\tau_j)\to0$ and $\varepsilon_j/\tau_j\to\lambda\in(0,+\infty]$ admit a subsequence {\rm(}still denoted by $(\varepsilon_j,\tau_j)${\rm)} and a limit function $u\colon[0,T]\mapsto X$ such that 
\begin{equation}\label{limite}
\tilde{u}_{\epsilon_j,\tau_j}(t)\to u(t)
\end{equation}
 for all $t \in [0, T]$. Moreover, $u$ satisfies the following properties:
\begin{enumerate}
\item[\rm(i)] $u$ is regulated;
\item[\rm(ii)] it holds
\begin{equation}\label{stazionario}
\nabla_x F(t,u_+(t))=\nabla_xF(t, u_-(t))=0 \quad \text{in $X$ for every $t\in [0,T]$};
\end{equation}
\item[\rm(iii)] denoting with $J$ the jump set of $u$, it holds $F(t, u_{+}(t))< F(t, u_{-}(t))$ for all $t\in J$.
\end{enumerate}
\label{compactnessthm}
\end{thm}

\proof

As we proved in the previous section, namely with \eqref{compactfunc}, the sequence $(\tilde{u}_{\epsilon,\tau}(t))_{\epsilon, \tau}$ is compact at every $t\in[0,T]$. {We denote by%Defining the positive measures
%\begin{equation*}
%\nu_{\varepsilon,\tau}:=\left(\frac{1}{4\varepsilon}\|\nabla_xF_{\tau}(\cdot,\tilde{u}_{\epsilon,\tau}(\cdot))\|^2+\frac{\varepsilon}{4}\|\dot{\bar{u}}_{\epsilon,\tau}(\cdot)\|^2\right)\,\mathcal{L}^1, 
%\end{equation*} 
\begin{equation*}
\nu_{\varepsilon,\tau}:=\left(\frac{1}{4\varepsilon}\|\nabla_xF_{\tau}(\cdot,\tilde{u}_{\epsilon,\tau}(\cdot))\|^2+\frac{\varepsilon}{4}\|\dot{\bar{u}}_{\epsilon,\tau}(\cdot)\|^2\right)\,\mathcal{L}^1
\end{equation*} 
the positive measures absolutely continuous with respect to the 1-dimensional Lebesgue measure $\mathcal{L}^1$ with density
$\frac{1}{4\varepsilon}\|\nabla_xF_{\tau}(\cdot,\tilde{u}_{\epsilon,\tau}(\cdot))\|^2+\frac{\varepsilon}{4}\|\dot{\bar{u}}_{\epsilon,\tau}(\cdot)\|^2$.}
Then from \eqref{enestimateimproved} we get that $\nu_{\varepsilon,\tau}$ are equibounded. Therefore, there exist sequences $(\varepsilon_j,\tau_j)\to0$ %with $\varepsilon_j/\tau_j\to+\infty$ 
as in the statement such that $\nu_{\varepsilon_j,\tau_j}$ converge weakly* to a positive finite measure $\nu$ on $[0,T]$. Then, the set of atoms $J_\nu$ of $\nu$ is at most countable {(see, e.g., \cite[Example~1.5(b)]{AFP})}. It also follows from \eqref{enestimateimproved} and the previous bounds that
\begin{align}\label{zero}
\nabla_xF_{\tau_j}(t,\tilde{u}_{\epsilon_j,\tau_j}(t))\to 0,\,\text{ as $j\to+\infty$, }\quad\mbox{ for a.e. }t\in [0,T]\,.
\end{align}

We may now fix a countable dense subset $I \subset [0,T]$ with the property that $I \supset J_\nu \cup\{0\}$, and define for all $t\in I$ the pointwise limit $u(t)$ of $\tilde{u}_{\epsilon_j,\tau_j}(t)$ (along a time independent subsequence) via a {diagonal} argument. We now want to show that $u$ can be {extended} to $[0, T]\setminus I$ in a unique way. Indeed, if $t \in [0,T]\setminus I$,  it holds in particular $t \notin J_\nu$. Let now $(t_l)_l$ and $(s_l)_l$ be two distinct sequences of points in the set $I$, both converging to $t$, and $u_1$ and $u_2$ be limit points of $u(t_l)$ and $u(s_l)$, respectively. Assume by contradiction that $u_1\neq u_2$. With a {diagonal} procedure we can extract subsequences $\tilde{u}_{\epsilon_l,\tau_l}$ of $ \tilde{u}_{\epsilon_j,\tau_j}$ such that
\begin{equation*}
\tilde{u}_{\epsilon_l,\tau_l}(t_l)\to u_1\quad\mbox{ and }\quad\tilde{u}_{\epsilon_l,\tau_l}(s_l)\to u_2\,.
\end{equation*}
Taking the corresponding nodes in $\Pi_{\tau_{l}}$ according to definition \eqref{node}, we get two sequences $\hat{t}_{\tau_{l}}$ and $\hat{s}_{\tau_{l}}$ both converging to $t$ such that, in particular,
\begin{equation}\label{diag}
\bar{u}_{\epsilon_l,\tau_l}(\hat{t}_{\tau_{l}})=\tilde{u}_{\epsilon_l,\tau_l}(\hat{t}_{\tau_{l}})\to u_1\quad\mbox{ and }\quad\bar{u}_{\epsilon_l,\tau_l}(\hat{s}_{\tau_{l}})=\tilde{u}_{\epsilon_l,\tau_l}(\hat{s}_{\tau_{l}})\to u_2\,.
\end{equation}
Denoting with $Q_l$ the nontrivial interval between  the two possibilities $[{\hat{s}_{\tau_{l}}}, \hat{t}_{\tau_{l}}]$ and $[{\hat{t}_{\tau_{l}}}, \hat{s}_{\tau_{l}}]$, we prove the following\\
{
\noindent
{\bf Claim:} there exist a subinterval $[a_l,b_l]$ of $Q_l$ and constants $\alpha,\delta>0$ independent of $l$ such that, for $l$ large enough, 
\begin{equation}
\|\bar{u}_{\epsilon_l,\tau_l}(b_l)-\bar{u}_{\epsilon_l,\tau_l}(a_l)\|>\delta\,\text{ and }\, \|\nabla_x F(s, \bar{u}_{\epsilon_l,\tau_l}(s))\|\geq\alpha,
\label{alphadelta}
\end{equation}
for every $s\in[a_l,b_l]$.}\\
{
For this, we notice that if $M>0$ is a uniform bound for $\bar{u}_{\epsilon_l,\tau_l}$, then by assumptions {\bf (F1)} and {\bf (F3)} the set $B_M\cap \mathcal{C}(t)$ is finite. Thus, there exists $\bar{\eta}=\bar{\eta}(t,M, u_1,u_2)$ such that, for every $0<\eta\leq \bar{\eta}$, it holds
\begin{equation}
B_{2\eta}(v)\cap B_{2\eta}(w)=\emptyset, \quad \text{ for every }v,w\in\mathcal{S},\, v\neq w,
\label{dupalle}
\end{equation}
where $\mathcal{S}=\mathcal{S}(t,u_1,u_2,M):=(B_M\cap \mathcal{C}(t))\cup\{u_1,u_2\}$. Now, if we introduce the compact set $K_\eta$ defined by $K_\eta:=\overline{B_M\backslash \bigcup_{v\in \mathcal{S}}B_\eta(v)}$, we have that $\displaystyle\min_{u\in K_\eta}\|\nabla_x F(t,u)\|>0$ and, by the regularity assumption {\bf (F0)}, there exists $\gamma=\gamma(t,M, u_1,u_2)>0$ such that
\begin{equation}
m_\eta:=\min_{u\in K_\eta, r\in[t-\gamma,t+\gamma]} \|\nabla_x F(r,u)\|>0.
\label{mineta}
\end{equation}
Since $\hat{s}_{\tau_{l}}\to t$ and $\hat{t}_{\tau_{l}}\to t$, for every $l$ sufficiently large we have that $Q_l\subset [t-\gamma,t+\gamma]$. Moreover, since $\bar{u}_{\epsilon_l,\tau_l}(\hat{s}_{\tau_{l}})\to u_1$ and $\bar{u}_{\epsilon_l,\tau_l}(\hat{t}_{\tau_{l}})\to u_2$, and from the definition of $K_\eta$ we also get that the set
\begin{equation*}
\mathcal{T}_l:=\{r\in Q_l:\,\bar{u}_{\epsilon_l,\tau_l}(r)\in K_\eta\}
\end{equation*}
is nonempty, for $l$ sufficiently large, and that there exist $a_l,b_l\in \mathcal{T}_l$, with $a_l\leq b_l$, such that $\|\bar{u}_{\epsilon_l,\tau_l}(a_l)-u_1\|=\eta$ and $\|\bar{u}_{\epsilon_l,\tau_l}(b_l)-u_2\|=\eta$. Thus, by \eqref{mineta} we get
\begin{equation*}
\|\nabla_x F(s, \bar{u}_{\epsilon_l,\tau_l}(s))\|\geq m_\eta=:\alpha,\,\text{ for every }s\in[a_l,b_l],
\end{equation*}
and
\begin{equation*}
\begin{split}
\|\bar{u}_{\epsilon_l,\tau_l}(b_l)-\bar{u}_{\epsilon_l,\tau_l}(a_l)\|&\geq \|u_1-u_2\|-(\|\bar{u}_{\epsilon_l,\tau_l}(a_l)-u_1\|+\|\bar{u}_{\epsilon_l,\tau_l}(b_l)-u_2\|)\\
&\geq \min_{v,w\in\mathcal{S}}(\|v-w\|-2\eta)=:\delta,
\end{split}
\end{equation*}
where $\delta$ is strictly positive by \eqref{mineta} and \eqref{dupalle}. This concludes the proof of Claim.} 

Coming back to the main proof, we have by the convergence of the positive measures $\nu_{\epsilon_{l},\tau_{l}}$ to $\nu$ and the Cauchy inequality
\begin{align*}
\nu([t-\eta, t+\eta])\ge \limsup_{l\to+\infty}\int_{Q_l} \left(\frac{\epsilon_l}4\|\dot {\bar{u}}_{\epsilon_l,\tau_l}(r)\|^2+\frac{1}{4\epsilon_l}\|\nabla_x F_{\tau_{l}}(r, \tilde{u}_{\epsilon_l,\tau_l}(r))\|^2\right)\mathrm{d}r \\
\ge \limsup_{l\to+\infty}  \frac{1}{2}\int_{Q_l}\|\dot {\bar{u}}_{\epsilon_l,\tau_l}(r)\|\|\nabla_x F_{\tau_{l}}(r, \tilde{u}_{\epsilon_l,\tau_l}(r))\|\mathrm{d}r 
\end{align*}
for any $\eta>0$. Letting $\eta \to 0$, since $t \notin J_\nu$ we deduce that
\begin{equation}\label{absurd}
0=\lim_{l\to+\infty} \int_{Q_l}\|\dot {\bar{u}}_{\epsilon_l,\tau_l}(r)\|\|\nabla_x F_{\tau_{l}}(r, \tilde{u}_{\epsilon_l,\tau_l}(r))\|\mathrm{d}r\,.
\end{equation}

{Now, let $[a_l,b_l]\subset Q_l$ and $\delta, \alpha>0$ be as in Claim. In particular, by virtue of \eqref{alphadelta}, the assumptions of Lemma~\ref{lemmafogli} hold. Thus, there exists a suitable $\beta>0$ such that $\|\nabla_x F_{\tau_{l}}(s, \tilde{u}_{\epsilon_l,\tau_l}(s))\|\geq \beta$ for every $s\in[a_l,b_l]$.}
%For $\beta$ as in Lemma~\ref{lemmafogli}, 
We then get 
\begin{equation}\label{decisiva}
\begin{split}
\int_{Q_l}\|\dot {\bar{u}}_{\epsilon_l,\tau_l}(r)\|\|\nabla_x F_{\tau_{l}}(r, \tilde{u}_{\epsilon_l,\tau_l}(r))\|\mathrm{d}r&\geq\int_{a_l}^{b_l}\|\dot {\bar{u}}_{\epsilon_l,\tau_l}(r)\|\|\nabla_x F_{\tau_{l}}(r, \tilde{u}_{\epsilon_l,\tau_l}(r))\|\mathrm{d}r\\
&\geq\beta\delta,
\end{split}
\end{equation}
and this gives a contradiction with \eqref{absurd}. Therefore it must be $u_1=u_2$. Setting $u(t)=u_1$ we can then extend $u$ in a unique way to a function defined on all $[0,T]$. 

The same argument as above can be now also used to show that
\begin{enumerate}
\item $\tilde{u}_{\epsilon_j,\tau_j}(t)$ converge to $u(t)$ also for $t\in [0,T]\backslash I$, and thus for every $t \in [0, T]$. 
\item $u$ is continuous at any $t\in [0,T]\backslash J_\nu$. Therefore, the jump set  $J$ of $u$ is contained in $J_\nu$, and it is at most countable. 
\end{enumerate}
Indeed, to prove the first claim, observe that for fixed $t\in [0,T]\backslash I$ one may now find a sequence $(t_l)_l$ of points in $I$ with $u(t_l)\to u(t)$ and take $(s_l)_l$ as being the sequence constantly equal to $t$. If (along a subsequence) $\tilde{u}_{\epsilon_j,\tau_j}(t)$ were not convergent to $u(t)$, we would get a contradiction by \eqref{diag}, \eqref{absurd}, and \eqref{decisiva}.
For proving the second claim, take $(t_l)_l$ as arbitrary sequence converging to $t \notin J_\nu$, and $(s_l)_l\equiv t$. Since now we know that $\tilde{u}_{\epsilon_j,\tau_j}(t)\to u(t)$, if $u(t_l)$ were not convergent to $u(t)$, we would get a contradiction by the fact that $t \notin J_\nu$ together with \eqref{diag}, \eqref{absurd}, and \eqref{decisiva}.

By pointwise convergence, \eqref{zero} and \eqref{(*)}, we also have $\nabla_x F(t, u(t))= 0$ for almost every $t\in [0,T]$. By continuity, equality holds then at every $t\in [0,T]\setminus J$, that is \eqref{staz} holds.

We now prove (i). To this aim we fix two sequences $(t_l)_l$ and $(s_l)_l$ with $t_l$, $s_l\searrow t$.  Let $u_1$ and $u_2$ be the limits of $u(t_l)$ and $u(s_l)$, respectively. Again we can  extract a subsequence $\tilde{u}_{\epsilon_l,\tau_l}$ such that
\[
\tilde{u}_{\epsilon_l,\tau_l}(t_l)\to u_1\mbox{ and }\tilde{u}_{\epsilon_l,\tau_l}(s_l)\to u_2\,,
\] 
From {\bf (F0)} and \eqref{(*)}, we have
\begin{align}\label{stim+}
\lim_{l\to +\infty}  F_{\tau_l}(t_l, \tilde{u}_{\epsilon_l,\tau_l}(t_l))- F(t_l, u(t_l))+ F_{\tau_l}(s_l, \tilde{u}_{\epsilon_l,\tau_l}(s_l))- F(s_l, u(s_l))=0\,.
\end{align}
We note that, setting $\hat{s}_{\tau_l}$ and $\hat{t}_{\tau_l}$ as in \eqref{node}, with \eqref{(*)} we get
\begin{equation}
\Bigl|F_{\tau_l}(s_l, \tilde{u}_{\epsilon_l,\tau_l}(s_l))-F_{\tau_l}(\hat{s}_{\tau_l}, \tilde{u}_{\epsilon_l,\tau_l}(\hat{s}_{\tau_l}))+F_{\tau_l}(t_l, \tilde{u}_{\epsilon_l,\tau_l}(t_l))-F_{\tau_l}(\hat{t}_{\tau_l}, \tilde{u}_{\epsilon_l,\tau_l}(\hat{t}_{\tau_l}))\Bigr|\to0,
\label{(3.23bis)}
\end{equation}
as $l\to+\infty$. 
Observe  also that, by Proposition~\ref{equi-boundvariation}, $F(t, u(t))$ has bounded variation and in particular admits a right limit at every $t$. 
Therefore we have
\begin{equation}\label{furtherstep}
\lim_{l\to +\infty}F(t_l, u(t_l))- F(s_l, u(s_l))=0\,.
\end{equation}
Up to further extraction, it holds either $\hat{t}_{\tau_{l}}\le \hat{s}_{\tau_{l}}$ or $\hat{s}_{\tau_{l}} \le \hat{t}_{\tau_{l}}$ for all $l$. Assuming this last one is the case, it follows, also using \eqref{enestimateimproved}, \eqref{stim+}, \eqref{(3.23bis)}, \eqref{furtherstep}, {\bf (F2)} and the uniform bound on $\tilde{u}_{\epsilon,\tau}$, that
\begin{align*}
&
\limsup_{l\to +\infty} \frac{1}{2}\int_{\hat{s}_{\tau_l}}^{\hat{t}_{\tau_l}}\|\dot {\bar{u}}_{\epsilon_l,\tau_l}(r)\|\|\nabla_x F_{\tau_l}(r, \tilde{u}_{\epsilon_l,\tau_l}(r))\|\mathrm{d}r\le\\
&\limsup_{l\to+\infty}\int_{\hat{s}_{\tau_l}}^{\hat{t}_{\tau_j}} \left(\frac{\epsilon_l}4\|\dot {\bar{u}}_{\epsilon_l,\tau_l}(r)\|^2+\frac{1}{4\epsilon_l}\|\nabla_x F_{\tau_l}(r, \tilde{u}_{\epsilon_l,\tau_l}(r))\|^2\right)\mathrm{d}r\le
\\
&\limsup_{l\to +\infty}  -F_{\tau_l}(\hat{t}_{\tau_l}, \tilde{u}_{\epsilon_l,\tau_l}(\hat{t}_{\tau_l}))+ F_{\tau_l}(\hat{s}_{\tau_l}, \tilde{u}_{\epsilon_l,\tau_l}(\hat{s}_{\tau_l}))+\int_{\hat{s}_{\tau_l}}^{\hat{t}_{\tau_l}}\partial_r F(r,\tilde{u}_{\epsilon_l,\tau_l}(r))\,dr\le\\
&\limsup_{l\to +\infty}  -F_{\tau_l}(t_l, \tilde{u}_{\epsilon_l,\tau_l}(t_l))+ F_{\tau_l}(s_l, \tilde{u}_{\epsilon_l,\tau_l}(s_l))+\int_{\hat{s}_{\tau_l}}^{\hat{t}_{\tau_l}}\partial_r F(r,\tilde{u}_{\epsilon_l,\tau_l}(r))\,dr\le\\
&\lim_{l\to +\infty}  -F(t_l, u(t_l))+ F(s_l, u(s_l))+\int_{\hat{s}_{\tau_l}}^{\hat{t}_{\tau_l}}\partial_r F(r,\tilde{u}_{\epsilon_l,\tau_l}(r))\,dr=0\,.
\end{align*}
With this, arguing as in \eqref{decisiva} we get $u_1=u_2$. This proves the existence of $u_+(t)$ at every $t$. The existence of $u_-(t)$ can be proved along the same lines. Once (i) is established, \eqref{stazionario} follows immediately from \eqref{staz} and {\bf(F0)}.

To conclude the proof, observe that if $u_{+}(t)\neq u_{-}(t)$, we may fix as in \eqref{diag} two sequences $\hat{t}_{\tau_l} > \hat{s}_{\tau_l}$, both converging to $t$, with
\begin{equation*}
\bar{u}_{\epsilon_l,\tau_l}(\hat{t}_{\tau_l})=\tilde{u}_{\epsilon_l,\tau_l}(\hat{t}_{\tau_l})\to u_{+}(t)\quad\mbox{ and }\quad\bar{u}_{\epsilon_l,\tau_l}(\hat{s}_{\tau_l})=\tilde{u}_{\epsilon_l,\tau_l}(\hat{s}_{\tau_l})\to u_{-}(t)\,.
\end{equation*}
From the above convergences, {\bf(F0)}, \eqref{enestimateimproved} and  the Cauchy-Schwarz inequality we then deduce
\[
F(t, u_{-}(t))- F(t, u_{+}(t))\ge \limsup_{l\to +\infty} \frac{1}{2}\int_{{\hat{s}_{\tau_l}}}^{{\hat{t}_{\tau_l}}}\|\dot {\bar{u}}_{\epsilon_l,\tau_l}(r)\|\|\nabla_x F_{\tau_l}(r, \tilde{u}_{\epsilon_l,\tau_l}(r))\|\mathrm{d}r\,.
\]
The same argument leading to \eqref{decisiva} shows the existence of two positive constants $\beta, \delta >0$ with
\[
F(t, u_{-}(t))- F(t, u_{+}(t))\ge \frac12 \beta \delta >0\,,
\]
thus giving (iii) and concluding the proof.
\endproof

We now deduce the two announced results about convergence of the piecewise affine interpolations and the stability condition when the ratio $\frac{\epsilon}{\tau}$ converges to a finite limit.

\begin{cor}\label{samelimit}
Under the assumptions of {\rm Theorem~\ref{compactnessthm}},  let $u(t)$ be the function given by \eqref{limite}, with at most countable jump set $J$. 
Then, the piecewise affine interpolations $\bar{u}_{\epsilon_j, \tau_j}$ defined in \eqref{piecewiseaffine} satisfy
\[
\bar{u}_{\epsilon_j, \tau_j}(t)\to u(t)
\]
for all $t \in [0, T] \setminus J$.
\end{cor}

\proof
By \eqref{confrontz}, \eqref{parallel} and the boundedness of $\frac{\epsilon}{\tau}$ we have for all $t$:
\[
\|\bar{u}_{\epsilon_j, \tau_j}(t)-\tilde{u}_{\epsilon_j, \tau_j}(t)\|\le C \|\nabla_x F_{\tau_j}(t, \tilde{u}_{\epsilon_j, \tau_j}(t))\|\,.
\]
If $t\notin J$ the right-hand side is vanishing by {\bf(F0)}, \eqref{(*)},  \eqref{staz}, and \eqref{limite}, proving the statement.
\endproof

\begin{prop}\label{improve}
Under the assumptions of {\rm Theorem~\ref{compactnessthm}},  let $u(t)$ be the regulated function given by \eqref{limite}.
Suppose additionally that
\begin{equation*}
0<\lim_{j\to +\infty}\frac{\epsilon_j}{\tau_j}=\lambda<+\infty
\end{equation*}
Then there hold
\begin{equation*}
F(t, u_+(t))\le F(t, v)+\frac{\lambda}2\|v-u_+(t)\|^2
\end{equation*}
and
\[
F(t, u_-(t))\le F(t, v)+\frac{\lambda}2\|v-u_-(t)\|^2
\]
 for all $t\in [0,T]$ and all $v\in X$.
\end{prop}
\proof
For all $t$ and all $v$ we have, using  \eqref{scheme} and \eqref{piecewiseconstant}
\begin{equation}\label{rewr}
F_ {\tau_j}(t, \tilde{u}_{\epsilon_j, \tau_j}(t))\le F_ {\tau_j}(t,v)+\frac{\epsilon_j}{2\tau_j}\|v-\tilde{u}_{\epsilon_j, \tau_j}(t-\tau_j)\|^2\,.
\end{equation}
By \eqref{piecewiseconstant}, \eqref{parallel} and the boundedness of $\frac{\epsilon}{\tau}$ we have for all $t$:
\[
\| \tilde{u}_{\epsilon_j, \tau_j}(t)-\tilde{u}_{\epsilon_j, \tau_j}(t-\tau_j)\|
\le C \|\nabla_x F_{\tau_j}(t, \tilde{u}_{\epsilon_j, \tau_j}(t))\|\,.
\]
If $t\notin J$ the right-hand side is vanishing by {\bf(F0)}, \eqref{(*)},  \eqref{staz}, and \eqref{limite}. With this, taking the limit in \eqref{rewr} we have, again using {\bf(F0)}, \eqref{(*)}, and \eqref{limite} that
\[
F(t, u(t))\le F(t, v)+\frac{\lambda}2\|v-u(t)\|^2
\]
for all $t\in [0,T]\setminus J$. To complete the proof, it simply suffices to take one-sided limits when $t \in J$ and use  {\bf(F0)}.
\endproof

\section{Energy balance}\label{energybalance}

\subsection{General results}\label{generals}
A first goal that we want to achieve, independently of the limit ratio $\frac{\epsilon}{\tau}$, is showing  that the limit evolution $u(t)$ provided by Theorem~\ref{compactnessthm} satisfies a balance between the stored energy and  the power spent  along the evolution in an interval of time $[s, t]\subset [0, T]$, up to a (positive) dissipation cost which is concentrated on the jump set of $u$, or equivalently (due to (iii) in Theorem~\ref{compactnessthm}) on the jump set of the energy $t\to F(t, u(t))$.

Later on we will characterise the dissipation cost in dependence of the limit of the ratio $\frac{\epsilon}{\tau}$: when this limit is $+\infty$ (Section \ref{casoinfinito}), we will retrieve the notion of {\it Balanced Viscosity solution} introduced in the time-continuous setting in \cite{Ago-Rossi}, while a different transition cost will appear in Section \ref{casofinito} whenever $\frac{\epsilon}{\tau}$ tends to a {\it finite  strictly positive} value.

Still for this subsection we do not need to distinguish between the two cases, while a crucial role will be played by the additional assumptions {\bf(F4)} and {\bf(F5)} on the energy. 
The proof strategy aims at showing that, as a sole consequence of the stability condition, the sum between the energy dissipated at jumps and the gap in the stored energy at a time $t$ with respect to the initial one must be greater or equal than the power spent in $[0, t]$. The opposite inequality is instead a direct consequence of Proposition \ref{equi-boundvariation}, which stems out of the minimizing movement scheme via compactness.

Adapting to our setting some arguments in \cite[Section 7]{MinSav},  we show at first that the right continuous representative of the nonincreasing function
\begin{align}\label{decr}
{f(t):=F(t,u(t))-\int_0^t\partial_rF(r,u(r))\,dr,}
\end{align}
considered in Proposition \ref{equi-boundvariation}, has a positive right derivative at every point. This corresponds, roughly speaking, to be nondecreasing (and thus, constant) up to the energy dissipated at jumps.

\begin{prop}
Let $f:[0,T]\to\R$ be defined as in \eqref{decr}.
Then, for any $t\in[0,T]$, the Dini lower right derivative of the right-continuous representative $f_+$ at $t$ is non-negative; i.e., it holds
\begin{equation}
D_+{f_+}(t):=\mathop{\lim\inf}_{h\searrow0}\frac{f_+(t+h)-f_+(t)}{h}\geq0.
\label{liminfpos}
\end{equation}
\label{proposition}
\end{prop}

\proof
We note that it holds
\begin{equation*}
f_+(t)=F(t,u_+(t))-\int_0^t\partial_rF(r,u(r))\,dr.
\end{equation*}
We first prove that
\begin{equation}
\mathop{\lim\inf}_{h\searrow0}\frac{f_+(t+h)-f_+(t)}{h}=\mathop{\lim\inf}_{h\searrow0}\frac{F(t,u_+(t+h))-F(t,u_+(t))}{h}.
\label{4.7}
\end{equation}
Indeed, by a direct computation
\begin{equation*}
\begin{split}
&|f_+(t+h)-f_+(t)-F(t,u_+(t+h))+F(t,u_+(t))|\\
&=\left|\int_{t}^{t+h}\Bigl[\partial_rF(r,u(r))-\partial_rF(r,u_+(t+h))\Bigr]\,dr\right|=o(h)
\end{split}
\end{equation*}
as $h\searrow0$, since $u_+$ is a right-continuous representative of $u$. Recall that \eqref{stazionario} gives
\begin{equation}\label{ct}
\nabla_xF(t+h,u_+(t+h))=0,\quad \text{ for every }\, t\in[0,T],\, h>0\,.
\end{equation}
Moreover, from assumption {\bf (F4)}, the fact that $u_+(t+h)\to u_+(t)$ as $h\to0$ and $u_+(t)\in \mathcal{C}(t)$ it follows that
\begin{equation*}
\displaystyle \mathop{\lim\inf}_{h\searrow0}\frac{F(t,u_+(t+h))-F(t,u_+(t))}{\|\nabla_xF(t,u_+(t+h))\|}\geq0.
\end{equation*}
We have
\begin{equation*}
\begin{split}
&\displaystyle \mathop{\lim\inf}_{h\searrow0}\frac{F(t,u_+(t+h))-F(t,u_+(t))}{h}\\
&\displaystyle \mathop{\lim\inf}_{h\searrow0}\frac{F(t,u_+(t+h))-F(t,u_+(t))}{\|\nabla_xF(t,u_+(t+h))\|}\frac{\|\nabla_xF(t,u_+(t+h))\|}{h},
\end{split}
\end{equation*}
where the term $\frac{\|\nabla_xF(t,u_+(t+h))\|}{h}$ is positive and bounded by assumption {\bf (F5)} and \eqref{ct}. From this we easily get
\begin{equation*}
\displaystyle \mathop{\lim\inf}_{h\searrow0}\frac{F(t,u_+(t+h))-F(t,u_+(t))}{h}\geq0,
\end{equation*}
which combined with \eqref{4.7} finally gives the inequality \eqref{liminfpos}.
\endproof

In order to prove our main result, we need the following elementary lemma, whose proof is only reported for the reader's convenience.

\begin{lem}
Let $g:[a,b]\to\R$ be continuous and such that the \emph{Dini upper right derivative of $g$ at $t$},
\begin{equation}
D^+g(t):=\mathop{\lim\sup}_{h\searrow0}\frac{g(t+h)-g(t)}{h}\geq0,\quad \forall t\in(a,b). 
\label{upperderivative}
\end{equation}
Then $g$ is nondecreasing on $(a,b)$.
\label{lem}
\end{lem}

\proof
It suffices to prove the result under the stronger assumption that $D^+g(t)>0$. Indeed, the general case follows by first replacing $g$ with $g_\eta(t):=g(t)+\eta t$ for some $\eta>0$, since $D^+g_\eta(t)=D^+g(t)+\eta>0$ implies $g_\eta$ increasing and then $g$ increasing by the arbitrariness of $\eta$. 

Now we argue by contradiction and assume the existence of $t_1,t_2\in(a,b)$, $t_1<t_2$, with $g(t_1)>g(t_2)$. Then, there exist $w$ such that $g(t_1)>w>g(t_2)$ and some points $s\in[t_1,t_2]$ such that $g(s)>w$. Denote by $\bar{s}:=\sup\{s\in[t_1,t_2]:\, g(s)>w\}$. We have $\bar{s}\in(t_1,t_2)$ and, by the continuity of $g$, $g(\bar{s})=w$. Therefore, for every $t\in(\bar{s},t_2)$
\begin{equation*}
\frac{g(t)-g(\bar{s})}{t-\bar{s}}<0,
\end{equation*}
from which follows $D^+g(\bar{s})\leq0$, a contradiction.
\endproof

\begin{oss}
The statement of Lemma~\ref{lem} holds true if we replace (\ref{upperderivative}) by an analogous assumption on the Dini lower right derivative of $g$ at $t$:
\begin{equation*}
D_+g(t):=\mathop{\lim\inf}_{h\searrow0}\frac{g(t+h)-g(t)}{h}\geq0,\quad \forall t\in(a,b), 
\end{equation*}
since the latter implies \eqref{upperderivative}.
\end{oss}

The following theoretical result provides the energy balance equality \eqref{measure} for our energies $F(t,u)$.

\begin{thm}\label{genbalance}
There exists a positive atomic measure $\mu$, with ${\rm supp}(\mu)=J$, $J$ being the jump set of $u$, such that
\begin{equation}
F(t,u_+(t))+\mu([s,t])=F(s,u_-(s))+\int_s^t \partial_r F(r,u(r))\,dr,
\label{measure}
\end{equation} 
for all $0\leq s\leq t\leq T$.
\end{thm}

\proof

Let $f: [0,T]\to\mathbb{R}$ be defined as in \eqref{decr}. We note that $f$ is  nonincreasing by Proposition~\ref{equi-boundvariation}. 
If we define $\mu:=-Df$, since $f$ is nonincreasing we have that $\mu$ is a positive measure satisfying $\mu([s,t])+Df([s,t])=0$, for any $s\leq t$, which corresponds to \eqref{measure} by virtue of \eqref{teoremafond}. We are left to show that ${\rm supp}(\mu)=J$. In order to do that, we define
\begin{equation*}
f^J(t):=\sum_{s\in[0,t]} (f_+(s)-f_-(s)),
\end{equation*}
which is the right-continuous jump function of $f$. We note that the set of discontinuities of $f^J$ coincides with $J$ and $Df^J=(Df)^J$, the latter being the jump part of measure $Df$. Moreover, $f^J$ is nonincreasing, so that $\mu^J=-Df^J$ is positive. It holds also $\mu\geq \mu^J$, since $\mu$ is positive. Now $-f^J$ is nondecreasing, so that combining with Proposition~\ref{proposition} we get
\begin{equation*}
\mathop{\lim\inf}_{h\searrow0}\frac{(f_+-f^J)(t+h)-(f_+-f^J)(t)}{h}\geq0.
\end{equation*}
Since, by construction, $f_+-f^J$ is a continuous function,  by Lemma~\ref{lem} $f_+-f^J$ is nondecreasing. Therefore $f_+(t)-f^J(t)\geq f_+(0)-f^J(0)$ for all $t\in[0,T]$, or, equivalently, 
\begin{equation*}
F(t,u_+(t))+\mu^J([0,t])\geq F(0,u(0))+\int_0^t \partial_r F(r,u(r))\,dr,
\end{equation*}
where, by the usual convention, $u(0)=u_-(0)$. Comparing the latter estimate with \eqref{measure} we finally get
\begin{equation*}
\mu^J([0,t])\geq\mu([0,t])\geq\mu^J([0,t]),\quad \forall t,
\end{equation*}
which gives $\mu^J=\mu$, thus concluding the proof.
\endproof

\begin{oss}
We note that, by construction, it holds
\begin{equation}\label{mu}
F(t,u_-(t))-F(t,u_+(t))=\mu(\{t\})>0,\quad \forall t\in J.
\end{equation}
\end{oss}

\subsection{The regime $\displaystyle\lim_{j\to+\infty}\frac{\epsilon_j}{\tau_j}=+\infty$.}\label{casoinfinito}

In this subsection, we will assume that
\begin{equation}\label{raf}
\lim_{j\to+\infty}\frac{\epsilon_j}{\tau_j}=+\infty
\end{equation}
and show, that under this assumption, the dissipation cost at each discontinuity time coincides with the one introduced in \cite{Ago-Rossi}. Therefore, the limit evolution $u(t)$ is a Balanced Viscosity solution in the sense of \cite[Definition 2]{Ago-Rossi}.

We introduce the \emph{cost function} $c_t:X\times X\to[0,+\infty)$, $t\in[0,T]$, defined by 
\begin{equation}
c_t(u_1,u_2):=\inf\left\{\int_0^1\|\dot{\theta}(s)\|\|\nabla_x F(t, \theta(s))\|\,ds:\,\theta\in{\rm AC}([0,1];X),\,\theta(0)=u_1, \theta(1)=u_2\right\},
\label{costfun}
\end{equation}
where ${\rm AC}([0,1];X)$ denotes the space of absolutely continuous functions defined on $[0,1]$ and taking values in $X$. {Notice that in \cite{Ago-Rossi} the cost is actually defined on a slightly larger class of continuous and piecewise absolutely continuous curves, but this clearly does not change the infimum. We remark however that existence of a minimum is in general only achieved in this larger class. For our purposes, this is nevertheless not relevant, since our proofs do not require the infimum above to be a minimum}.
We recall without proof \cite[Theorem~2.4]{Ago-Rossi}, which collects the main properties of the cost function $c_t$.

\begin{thm}[Agostiniani-Rossi~\cite{Ago-Rossi}] Assume {\bf (F0)-(F3)}. Then, for every $t\in[0,T]$ and $u_1,u_2\in X$, the function $c_t(u_1, u_2)$ defined by \eqref{costfun} satisfies:
\begin{enumerate}
\item[\rm(i)] $c_t(u_1,u_2)=0$ if and only if $u_1=u_2$;
\item[\rm(ii)] $c_t$ is symmetric;
\item[\rm(iii)] if $c_t(u_1,u_2)>0$, there exists an optimal curve $\theta$ attaining the inf in {\rm(\ref{costfun})};
\item[\rm(iv)] for every $u_3\in \mathcal{C}(t)$, the triangle inequality holds
\begin{equation*}
c_t(u_1,u_2)\leq c_t(u_1,u_3)+c_t(u_3,u_2);
\end{equation*} 
\item[\rm(v)] there holds 
\begin{equation*}
\begin{split}
c_t(u_1,u_2)\leq \inf &\Bigl\{\mathop{\lim\inf}_{j\to+\infty}\int_{t^1_j}^{t^2_j}\|\dot{\theta}_j(s)\|\|\nabla_x F(s, \theta_j(s))\|\,ds:\\
&\theta_j\in{\rm AC}([t^1_j,t^2_j];X),\, t^i_j\to t,\,\theta_j(t^i_j)\to u_i\, \text{ for }i=1,2\Bigr\};
\end{split}
\end{equation*}
\item[\rm(vi)] the following lower semicontinuity property holds
\begin{equation*}
(u_1^l,u_2^l)\to (u_1,u_2)\,\text{ as }l\to+\infty \Rightarrow \mathop{\lim\inf}_{l\to+\infty}c_t(u_1^l,u_2^l)\geq c_t(u_1,u_2).
\end{equation*}
\end{enumerate}
\label{thm2.4}
\end{thm} 

We start  with a technical Lemma, which essentially stems out of Proposition \ref{equi-boundvariation}. Its heuristical meaning is that the stored energy is asymptotically decreasing along the fast transition from $u_-(t)$ to $u_+(t)$ since the power spent is vanishing by the absolute continuity of the integral. This will eventually allow us, in the proof of \eqref{lowbounddiss}, to cut out the intervals where
\begin{enumerate}
\item
either $\bar{u}_{\epsilon,\tau}$ stays close to $\mathcal C(t)$, so that no good control on $\dot{\bar{u}}_{\epsilon,\tau}$ is available and Lemma \ref{lemma2} does not hold,
\item
or $\bar{u}_{\epsilon,\tau}$ makes a loop from and to the same point in $\mathcal C(t)$. This is non-optimal since it forces additional energy to be dissipated.
\end{enumerate}

In the lemma below, we make also use of the additional hypothesis \eqref{uniforme}, which is verified along the sequence $(\epsilon_j, \tau_j)$ when \eqref{raf} holds, as a consequence of Proposition~\ref{tend0} (ii).
\begin{lem}
Let $[a,b]\subseteq[0,T]$ and $\displaystyle A=\bigcup_{i\in I}(a^i,b^i)\subset[a,b]$, with $b^i\leq a^{i+1}$, $I\subset\mathbb{N}$, $\#I=m<\infty$. 
Assume that
\begin{equation}\label{uniforme}
\|\bar{u}_{\epsilon,\tau}(t)-\tilde{u}_{\epsilon,\tau}(t)\|\to0
\end{equation}
uniformly with respect to $t$ in $[0,T]$ as $(\epsilon, \tau)\to (0,0)$.
Then there exist a constant $C>0$, independent of $\epsilon$ and $\tau$, and a remainder $\zeta_{\epsilon, \tau}$ such that $\zeta_{\epsilon, \tau}\to0$ when $(\epsilon, \tau)\to (0,0)$ for which it holds 
\begin{equation}
\begin{split}
F_{\tau}(a,\bar{u}_{\epsilon,\tau}(a))-F_{\tau}(b,\bar{u}_{\epsilon,\tau}(b))&\geq \sum_{i\in I} \left(F_{\tau}(a^i,\bar{u}_{\epsilon,\tau}(a^i))-F_{\tau}(b^i,\bar{u}_{\epsilon,\tau}(b^i))\right)\\
&-C(b-a+\tau)-\zeta_{\epsilon, \tau}(m+2).
\end{split}
\label{stimalemma}
\end{equation}
\label{lemma1}
\end{lem}

\proof

We denote by $\hat{a}^i_{\tau}$ and $\hat{b}^i_{\tau}$ the nodes of partition $\Pi_{\tau}$ corresponding to $a^i$ and $b^i$, respectively, defined as in \eqref{node}. We have $\hat{a}^i_{\tau}\leq \hat{b}^i_{\tau}$ and, as already remarked in \eqref{identity}, $\bar{u}_{\epsilon,\tau}(\hat{a}^i_{\tau})=\tilde{u}_{\epsilon,\tau}({a}^i)$ and $\bar{u}_{\epsilon,\tau}(\hat{b}^i_{\tau})=\tilde{u}_{\epsilon,\tau}({b}^i)$. By  \eqref{uniforme} and {\bf (F0)}, there exists a vanishing remainder $\zeta_{\epsilon, \tau}$ such that
\begin{equation}\label{1step}
\begin{split}
\sum_{i\in I} &\left|F_{\tau}(a^i,\bar{u}_{\epsilon,\tau}(a^i))-F_{\tau}(b^i,\bar{u}_{\epsilon,\tau}(b^i))-\left(F_{\tau}(\hat{a}^i_{\tau},\bar{u}_{\epsilon,\tau}(\hat{a}^i_{\tau}))-F_{\tau}(\hat{b}^i_{\tau},\bar{u}_{\epsilon,\tau}(\hat{b}^i_{\tau}))\right)\right|\\
&\leq m\zeta_{\epsilon, \tau}\,.
\end{split}
\end{equation}
For each $i$, set $\mathcal{K}^i_{\tau}:=\{ k \in \mathcal{K}_\tau:\hat{ b}^{i+1}_{\tau}\le t^k <\hat{a}^{i+1}_{\tau}\}$. We have, also using \eqref{stimaa} 
\begin{equation}\label{2step}
\begin{array}{c}
\displaystyle
F_{\tau}(\hat{a}^{i+1}_{\tau},\bar{u}_{\epsilon,\tau}(\hat{a}^{i+1}_{\tau}))-F_{\tau}(\hat{b}^i_{\tau},\bar{u}_{\epsilon, \tau}(\hat{b}^i_{\tau}))=
\sum_{k\in \mathcal{K}^i_{\tau}} \left(F(t^{k+1}, u^{k+1})- F(t^{k}, u^{k})\right) \le \\[10 pt]
\displaystyle
\sum_{k\in \mathcal{K}^i_{\tau}} \left(\int^{t^{k+1}}_{t^{k}}\partial_r F(r, u^{k})\,\mathrm{d}r\right) \le C (\hat{a}^{i+1}_{\tau}-\hat{b}^i_{\tau})
\end{array}
\end{equation}
for all $i \in I$, where $C$ is a positive constant such that $|\partial_r F(r,u^k)|\leq C$ for all $k\in\mathcal{K}_{\tau}$.
By construction, it holds $\hat{b}^i_{\tau}\leq\hat{a}^{i+1}_{\tau}$ for any $i=1,\dots,m$, which gives
\begin{equation*}
\sum_{i\in I}(\hat{a}^{i+1}_{\tau}-\hat{b}^i_{\tau})\leq b-a+2\tau.
\end{equation*}
With this, \eqref{1step} and \eqref{2step} we get
\begin{equation}\label{3step}
\begin{split}
\sum_{i\in I}\left(F_{\tau}(a^i,\bar{u}_{\epsilon,\tau}(a^i))-F_{\tau}(b^i,\bar{u}_{\epsilon,\tau}(b^i))\right)&\leq F_{\tau}(\hat{a}^1_{\tau},\bar{u}_{\epsilon,\tau}(\hat{a}^1_{\tau}))-F_{\tau}(\hat{b}^m_{\tau},\bar{u}_{\epsilon,\tau}(\hat{a}^m_{\tau}))\\
&+C(b-a+2\tau)+m\zeta_{\epsilon, \tau}.
\end{split}
\end{equation}

A further application of the argument in \eqref{2step} leads to
\begin{equation}\label{4step}
F_{\tau_j}(\hat{a}^1_{\tau},\bar{u}_{\epsilon,\tau}(\hat{a}^1_{\tau}))-F_{\tau}(\hat{a}_{\tau},\bar{u}_{\epsilon,\tau}(\hat{a}_{\tau}))\leq C(b-a+2\tau),
\end{equation}
and
\begin{equation}\label{5step}
F_{\tau_j}(\hat{b}_{\tau},\bar{u}_{\epsilon,\tau}(\hat{b}_{\tau}))-F_{\tau}(\hat{b}^m_{\tau},\bar{u}_{\epsilon,\tau}(\hat{b}^m_{\tau}))\leq C(b-a+2\tau).
\end{equation}
Finally, again by \eqref{uniforme} and {\bf (F0)}
\begin{equation*}
\Bigl|F_{\tau}(a,\bar{u}_{\epsilon,\tau}(a)-F_{\tau}(\hat{a}_{\tau},\bar{u}_{\epsilon,\tau}(\hat{a}_{\tau}))-F_{\tau}(b,\bar{u}_{\epsilon,\tau}(b)+F_{\tau}(\hat{b}_{\tau},\bar{u}_{\epsilon,\tau}(\hat{b}_{\tau}))\Bigr|\leq 2\zeta_{\epsilon,\tau},
\end{equation*}
With this,  \eqref{stimalemma} follows from \eqref{3step}, \eqref{4step}, and \eqref{5step},  up to redefining $C$ and $\zeta_{\epsilon, \tau}$.
\endproof

 The following Lemma is a key result for showing that $\bar{u}_{\epsilon_j,\tau_j}$ optimizes the cost $c_t(\cdot, \cdot)$ along a transition.

\begin{lem}\label{lemma2}
Let $\eta>0$ be fixed and define
\begin{equation*}
A_j^\eta:=\Bigl\{t\in[0,T]:\, {\rm dist}(\bar{u}_{\epsilon_j,\tau_j}(t),\mathcal{C}(t))>\eta\Bigr\},
\end{equation*}
Assume that \eqref{raf} holds. Then
\begin{equation}
\lim_{j\to+\infty}\int_{A_j^\eta}\Bigl|\langle\nabla_x F_{\tau_j}(s, \bar{u}_{\epsilon_j,\tau_j}(s)),\dot{\bar{u}}_{\epsilon_j,\tau_j}(s)\rangle-\|\nabla_x F(s, \bar{u}_{\epsilon_j,\tau_j}(s))\|\|\dot{\bar{u}}_{\epsilon_j,\tau_j}(s)\|\Bigr|\,ds=0.
\label{thlemma2}
\end{equation}
\end{lem}

\proof
Since $\|\nabla_x F(\cdot,\cdot)\|$ is continuous, corresponding to $\eta$ there exists a positive constant $\alpha_\eta$ such that
\begin{equation*}
\|\nabla_x F(t,u)\|\geq\alpha_\eta,\quad \text{ for every }\, (t,u)\, \text{ such that }\, {\rm dist}(u,\mathcal{C}(t))\geq\eta.
\end{equation*}
By the uniform convergence (\ref{(*)}), the same bound holds also for $\|\nabla_x F_{\tau_j}(t,u)\|$, up to eventually change $\alpha_\eta$. Thus, we have
\begin{equation*}
\|\nabla_x F_{\tau_j}(t,\bar{u}_{\epsilon_j,\tau_j}(t))\|\geq\alpha_\eta,\quad \text{ for every }\, t\in A_j^\eta,
\end{equation*}
and, as a consequence of Lemma~\ref{lemmafogli}, correspondingly we can find $\beta_\eta$ such that
\begin{equation*}
\|\nabla_x F_{\tau_j}(t,\tilde{u}_{\epsilon_j,\tau_j}(t))\|\geq\beta_\eta,\quad \text{ for every }\, t\in A_j^\eta.
\end{equation*}
With this, \eqref{enestimateimproved}, and the Cauchy inequality, we have
\begin{equation}\label{l1bound}
\int_{A_j^\eta}\|\dot{\bar{u}}_{\epsilon_j,\tau_j}(t)\|\,dt\leq\frac{1}{\beta_\eta}\int_{A_j^\eta}\|\nabla_x F_{\tau_j}(t,\tilde{u}_{\epsilon_j,\tau_j}(t))\|\|\dot{\bar{u}}_{\epsilon_j,\tau_j}(t)\|\,dt\leq \frac{C}{\beta_\eta}\,.
\end{equation}
Now \eqref{raf} and Proposition~\ref{tend0} (ii) imply that \eqref{uniforme} holds. With the equicontinuity of $\|\nabla_x F_{\tau}(\cdot,\cdot)\|$, which follows from {\bf(F0)}, we get
\begin{equation}\label{l8bound}
\|\nabla_x F_{\tau_j}(t,\tilde{u}_{\epsilon_j,\tau_j}(t))-\nabla_x F_{\tau_j}(t,\bar{u}_{\epsilon_j,\tau_j}(t))\|\to0
\end{equation}
uniformly in $[0,T]$ as $j\to+\infty$.  
Now, the  identity 
\begin{equation*}
-\langle\nabla_x F_{\tau_j}(t, \tilde{u}_{\epsilon_j,\tau_j}(t)),\dot{\bar{u}}_{\epsilon_j,\tau_j}(t)\rangle=\|\nabla_x F_{\tau_j}(t, \tilde{u}_{\epsilon_j,\tau_j}(t))\|\|\dot{\bar{u}}_{\epsilon_j,\tau_j}(t)\|
\end{equation*}
holds for a.e. $t\in[0,T]$ as a consequence of \eqref{parallel}. With this, \eqref{l1bound} and \eqref{l8bound} we get
\begin{equation*}
\lim_{j\to+\infty}\int_{A_j^\eta}\Bigl|\langle\nabla_x F_{\tau_j}(s, \bar{u}_{\epsilon_j,\tau_j}(s)),\dot{\bar{u}}_{\epsilon_j,\tau_j}(s)\rangle-\|\nabla_x F_{\tau_j}(s, \bar{u}_{\epsilon_j,\tau_j}(s))\|\|\dot{\bar{u}}_{\epsilon_j,\tau_j}(s)\|\Bigr|\,ds=0.
\end{equation*}
A further application of  the uniform convergence (\ref{(*)}) with the estimate \eqref{l1bound} yield finally \eqref{thlemma2}.
\endproof

We can now show that $c_t(u_+(t),u_-(t))$ is a lower bound for the dissipation $\mu(\{t\})$ at a jump point $t$.
\begin{prop}\label{taupiccolo}
Let $c_t$ be the cost function defined in \eqref{costfun}, $u_-(t)$ and $u_+(t)$ be the left and right limits, respectively, of the function $u$ of {\rm Theorem~\ref{compactnessthm}} at each point $t$.  Assume in addition that
\[
\lim_{j\to+\infty}\frac{\epsilon_j}{\tau_j}=+\infty\,.
\]
Then it holds
\begin{equation}
F(t,u_-(t))-F(t,u_+(t))\geq c_t(u_+(t),u_-(t)),\quad \forall t\in [0,T].
\label{lowbounddiss}
\end{equation}
\label{proposizione1}
\end{prop}

\proof
We restrict to the case of $t\in J$, since for any $t\in [0,T]\backslash J$ the equality holds as a consequence of (i) in Theorem \ref{thm2.4}.
First, it is not restrictive to assume that there are two sequences $t_j\searrow t$ and $s_j\nearrow t$ such that
\begin{equation}\label{fatto1}
\lim_{j\to +\infty}\|\bar{u}_{\epsilon_j,\tau_j}(s_j)- u_-(t)\|+ \|\bar{u}_{\epsilon_j,\tau_j}(t_j)- u_+(t)\|=0\,.
\end{equation}
Indeed, since $u$ is regulated and \eqref{limite} and \eqref{uniforme} hold, for two arbitrary sequences   $t_l\searrow t$ and $s_l\nearrow t$ a subsequence $\bar{u}_{\epsilon_{j_l},\tau_{j_l}}$ such that
\[
\lim_{j\to +\infty}\|\bar{u}_{\epsilon_{j_l},\tau_{j_l}}(s_l)- u_-(t)\|+ \|\bar{u}_{\epsilon_{j_l},\tau_{j_l}}(t_l)- u_+(t)\|=0
\]
can be constructed by a {diagonal} argument. {Maybe abusing a bit of notation}, we will avoid relabeling and assume \eqref{fatto1}. Notice that this latter implies
\begin{equation}\label{servira}
F_{\tau_j}(s_j,\bar{u}_{\epsilon_j,\tau_j}(s_j))-F_{\tau_j}(t_j,\bar{u}_{\epsilon_j,\tau_j}(t_j))\to F(t,u_-(t))-F(t,u_+(t)),
\end{equation}
as $j\to+\infty$.

For $M$ as in \eqref{compactfunc} and $\mathcal{C}(t)$ as in {\bf(F3)}, we set $\mathcal{C}_M(t):=\mathcal{C}(t)\cap B_M$. Notice that  {\bf(F3)} implies that $\mathcal{C}_M(t)$ is a finite set, and we denote by $N_t:=\#\mathcal{C}_M(t)$ its cardinality. We then define the strictly positive value $d$ as
\begin{equation*}
d=d_t:=\min\Bigl\{\|w-z\|:\, w,z\in \mathcal{C}_M(t),\,w\neq z\Bigr\}
\end{equation*}
and we fix $\eta>0$ %, %$\eta<\!\!<d$. 
with the property that
\begin{equation}
\eta<\frac{d}{2}.
\label{smalleta}
\end{equation}
With this,
\begin{equation*}
B_\eta(u^i)\cap B_\eta(u^j)=\emptyset,\quad \text{ for every } u^i,u^j\in\mathcal{C}_M(t),\quad i\neq j\,.
\end{equation*}
In particular, if for some $u\in B_M$, it holds
\begin{equation*}
\dist(u,\mathcal{C}_M(t))\leq\eta,
\end{equation*}
then there exists a unique $\hat{u}\in\mathcal{C}_M(t)$ such that
\begin{equation}
\dist(u,\mathcal{C}_M(t))=\|u-\hat{u}\|\leq\eta.
\label{uniqueproj}
\end{equation}

Define the open set $A_j^\eta\subset[0,T]$ as in Lemma~\ref{lemma2}. Our first aim is to prove the following\\
{\bf Claim:} for every $j\in\N$, there exists an open subset
\begin{equation*}
B_j^\eta:=\bigcup_{i\in I_j}(a_j^i,b_j^i)\subseteq A_j^\eta\cap(s_j,t_j),
\end{equation*}
such that $m_j:=\# I_j\leq N_t$, and a set of distinct critical points of $F(t,\cdot)$, say $\mathcal{U}:=\{u^0,u^1,\dots,u^{m_j}\}\subseteq \mathcal{C}_M(t)$, with $u^0=u_-(t)$ and $u^{m_j}=u_+(t)$, such that the following properties are satisfied:
\begin{enumerate}
\item[(1)] $\|\bar{u}_{\epsilon_j,\tau_j}(a_j^1)-u_-(t)\|=\eta$;
\item[(2)] $\|\bar{u}_{\epsilon_j,\tau_j}(b_j^{m_j})-u_+(t)\|=\eta$;
\item[(3a)] ${\rm dist}(\bar{u}_{\epsilon_j,\tau_j}(a_j^i),\mathcal{C}_M(t))=\|\bar{u}_{\epsilon_j,\tau_j}(a_j^i)-u^{i-1}\|=\eta$;
\item[(3b)] ${\rm dist}(\bar{u}_{\epsilon_j,\tau_j}(b_j^i),\mathcal{C}_M(t))=\|\bar{u}_{\epsilon_j,\tau_j}(b_j^i)-u^{i}\|=\eta$,\quad \text{ for every }\,$i\in I_j$;
\item[(4)] $\|\bar{u}_{\epsilon_j,\tau_j}(s)-u^{i-1}\|>\eta$\, for every $s>a_j^i$ and every $i\in I_j$.
%\item[(3)]$\|\bar{u}_{\epsilon_j,\tau_j}(b_j^i)-\bar{u}_{\epsilon_j,\tau_j}(a_j^{i+1})\|\leq2\eta$,\quad \text{ for every }\, $i=1,\dots,m_j-1$;
%\item[(4)]$\|\bar{u}_{\epsilon_j,\tau_j}(a_j^i)-\bar{u}_{\epsilon_j,\tau_j}(b_j^{i})\|\geq d-2\eta$,\quad \text{ for every }\,$i\in I_j$, 
%\item[(5)]${\rm dist}(\bar{u}_{\epsilon_j,\tau_j}(a_j^i),\mathcal{C}_M(t))=\eta$\, and\,${\rm dist}(\bar{u}_{\epsilon_j,\tau_j}(b_j^i),\mathcal{C}_M(t))=\eta$\, for all $i\in I_j$.
\end{enumerate}
Notice that (3a), (3b) and (4) together imply that $u^i\neq u^l$ for every $i,j$ with $i\neq l$. In order to prove the Claim, we will exhibit a construction of the set $B^\eta_j$ (and, correspondingly, of $\mathcal{U}$) arguing by induction on the number of intervals $m_j$. %It is worth to note that, by (4), each of $\bar{u}_{\epsilon_j,\tau_j}(a_j^i)$ and $\bar{u}_{\epsilon_j,\tau_j}(a_j^{i+1})$ will be close to a different point of $\mathcal{C}_M(t)$, so that the iterative procedure will stop after at most $N_t$ steps. Thus, $m_j\leq N_t$.
\\
{\bf Step 1.} Since $\|\bar{u}_{\epsilon_j,\tau_j}(s_j)-u_-(t)\|\to0$ as $j\to+\infty$ and
\begin{equation*}
\displaystyle\mathop{\lim\inf}_{j\to\infty}\|\bar{u}_{\epsilon_j,\tau_j}(t_j)-u_-(t)\|\geq\|u_+(t)-u_-(t)\|\geq d,
\end{equation*} 
it is well defined
\begin{equation*}
a_j^1:=\max\Bigl\{s\in[s_j,t_j]:\, \|\bar{u}_{\epsilon_j,\tau_j}(s)-u_-(t)\|\leq\eta\Bigr\}
\end{equation*}
and it satisfies $a_j^1<t_j$. Moreover
\begin{equation*}
{\rm dist}(\bar{u}_{\epsilon_j,\tau_j}(a_j^1),\mathcal{C}_M(t))=\|\bar{u}_{\epsilon_j,\tau_j}(a_j^1)-u_-(t)\|=\eta.
\end{equation*}
Now, setting
\begin{equation*}
b^1_j:=\min\Bigl\{s\in(a_j^1,t_j]:\, {\rm dist}(\bar{u}_{\epsilon_j,\tau_j}(s),\mathcal{C}_M(t))\leq\eta\Bigr\}
\end{equation*}
(note that the previous definition is well-posed), we have
\begin{equation*}
{\rm dist}(\bar{u}_{\epsilon_j,\tau_j}(b_j^1), \mathcal{C}_M(t))=\eta.
\end{equation*}
Moreover, with \eqref{smalleta}-\eqref{uniqueproj}, there exists a unique $u^1\in\mathcal{C}_M(t)$ such that
\begin{equation}
{\rm dist}(\bar{u}_{\epsilon_j,\tau_j}(b_j^1), \mathcal{C}_M(t))=\|\bar{u}_{\epsilon_j,\tau_j}(b_j^1)-u^1\|=\eta.
\label{(**)}
\end{equation}
On the other side, we note that
\begin{equation}
\|\bar{u}_{\epsilon_j,\tau_j}(b_j^1)-u_-(t)\|>\eta
\label{(**)2}
\end{equation}
since by the definition of $a_j^1$, for any $s>a_j^1$ it results $\|\bar{u}_{\epsilon_j,\tau_j}(s)-u_-(t)\|>\eta$. With \eqref{(**)} and \eqref{(**)2} we deduce that $u^1\neq u_-(t)$.
%With (\ref{(**)}) and the definition of $d$ we immediately get
%\begin{equation*}
%\|\bar{u}_{\epsilon_j,\tau_j}(b_j^1)-\bar{u}_{\epsilon_j,\tau_j}(a_j^1)\|\geq d-2\eta,
%\end{equation*}
%as desired. If $\|\bar{u}_{\epsilon_j,\tau_j}(b_j^1)-u_+(t)\|=\eta$, 
If $u^1=u_+(t)$, then we conclude and $m_j=1$. Otherwise, the construction goes on.
\\
{\bf Step 2:} Assume that $(a^l_j,b^l_j)\subset B_j^\eta$ and $u^l\in\mathcal{U}$ have been constructed for $1\leq l\leq i$ for some $i>1$, such that all the properties in the Claim are satisfied with the exception of (2). %Assume that $(a^i_j,b^i_j)\subset B_j^\eta$ is given for some $i>1$, so that $\#I_j=m_j-1$ for some integer $m_j\geq2$ and $\{u^0,u^1,\dots,u^{m_j-1}\}\subseteq\mathcal{U}$. 
Then we define
\begin{equation*}
a_j^{i+1}=\max\Bigl\{s\in [b_j^i,t_j]:\, \|\bar{u}_{\epsilon_j,\tau_j}(s)-{u}^i\|\leq \eta\Bigr\},
\end{equation*}
where ${u}^i\in\mathcal{C}_M(t)$ is such that
\begin{equation*}
{\rm dist}(\bar{u}_{\epsilon_j,\tau_j}(b_j^i),\mathcal{C}_M(t))=\|\bar{u}_{\epsilon_j,\tau_j}(b_j^i)-{u}^i\|.
\end{equation*}
It holds $a_j^{i+1}<t_j$, since by construction ${u}^i\neq u_+(t)$, and then
\begin{equation*}
\displaystyle\mathop{\lim\inf}_{j\to+\infty}\|\bar{u}_{\epsilon_j,\tau_j}(t_j)-{u}^i\|=\|u_+(t)-{u}^i\|\geq d.
\end{equation*}
Moreover, ${\rm dist}(\bar{u}_{\epsilon_j,\tau_j}(a_j^{i+1}),\mathcal{C}_M(t))=\|\bar{u}_{\epsilon_j,\tau_j}(a_j^{i+1})-{u}^i\|=\eta$. %and
%\begin{equation*}
%\|\bar{u}_{\epsilon_j,\tau_j}(b_j^{i})-\bar{u}_{\epsilon_j,\tau_j}(a_j^{i+1})\|\leq \|\bar{u}_{\epsilon_j,\tau_j}(b_j^{i})-\hat{u}_i\|+\|\hat{u}_i-\bar{u}_{\epsilon_j,\tau_j}(a_j^{i+1})\|=2\eta.
%\end{equation*}
Setting
\begin{equation*}
b_j^{i+1}=\min\Bigl\{s\in(a_j^{i+1},t_j]:\, {\rm dist}(\bar{u}_{\epsilon_j,\tau_j}(s), \mathcal{C}_M(t))\leq\eta\Bigr\}
\end{equation*}
%then conditions (4) and (5) are satisfied. 
as in Step 1 there exists $u^{i+1}\in\mathcal{C}_M(t)$ such that
\begin{equation}
{\rm dist}(\bar{u}_{\epsilon_j,\tau_j}(b_j^{i+1}), \mathcal{C}_M(t))=\|\bar{u}_{\epsilon_j,\tau_j}(b_j^{i+1})-u^{i+1}\|=\eta.
\label{etastim}
\end{equation}
We note also that the minimum in the previous equation is achieved since in a right neighborhood of $a_j^{i+1}$ we have
\begin{equation*}
 {\rm dist}(\bar{u}_{\epsilon_j,\tau_j}(s), \mathcal{C}_M(t))=\|\bar{u}_{\epsilon_j,\tau_j}(s)-{u}^i\|>\eta,
\end{equation*}
which, in particular, implies that
\begin{equation*}
\|\bar{u}_{\epsilon_j,\tau_j}(b_j^{i+1})-{u}^i\|>\eta.
\end{equation*}
The latter one combined with \eqref{etastim} gives $u^{i+1}\neq u^{l}$, for every $l\leq i$.

Thus, $(a_j^{i+1},b_j^{i+1})\subset B_j^\eta$, $m_j=i+1$ and $u^{i+1}\in\mathcal{U}$. With this, since $\mathcal{C}_M(t)$ has finite cardinality $N_t$ and recalling that $\bar{u}_{\epsilon_j,\tau_j}(t_j)\to u_+(t)$, in a finite number of steps $m_j\leq N_t$ we will get property (2), concluding the proof of the Claim. 

Let us go back to the main proof. Since $\#I_j\leq N_t$, up to passing to a subsequence, we may assume $\#I_j=m$ for any $j$, with $m$ independent of $j$. From Lemma~\ref{lemma1} and \eqref{servira} we get
\begin{equation*}
F(t,u_-(t))-F(t,u_+(t))\geq\displaystyle\sum_{i=1}^m\mathop{\lim\sup}_{j\to+\infty}\Bigl[F_{\tau_j}(a_j^i,\bar{u}_{\epsilon_j,\tau_j}(a_j^i))-F_{\tau_j}(b_j^i,\bar{u}_{\epsilon_j,\tau_j}(b_j^i))\Bigr].
\end{equation*}
The fundamental theorem of calculus coupled with Proposition~\ref{convergence1} gives now
\begin{equation*}
\begin{split}
F_{\tau_j}(a_j^i, \bar{u}_{\epsilon_j,\tau_j}(a_j^i))-F_{\tau_j}(b_j^i, \bar{u}_{\epsilon_j,\tau_j}(b_j^i))=&-\int_{a_j^i}^{b_j^i}\langle \nabla_x F_{\tau_j}(s,\bar{u}_{\epsilon_j,\tau_j}(s)),\dot{\bar{u}}_{\epsilon_j,\tau_j}(s)\rangle\,ds\\
&-\int_{a_j^i}^{b_j^i}\partial_r F(r,\tilde{u}_{\epsilon_j,\tau_j}(r))\,dr+o(1),
\end{split}
\end{equation*}
as $j\to+\infty$. Since $\partial_rF(r,\tilde{u}_{\epsilon_j,\tau_j}(r))$ is equi-bounded and $\displaystyle\sum_{i=1}^m(b_j^i-a_j^i)\leq t_j-s_j\to0$ as $j\to+\infty$, we get
\begin{equation}
F(t,u_-(t))-F(t,u_+(t))\geq \displaystyle\sum_{i=1}^m\mathop{\lim\sup}_{j\to+\infty}\left(-\int_{a_j^i}^{b_j^i}\langle \nabla_x F_{\tau_j}(s,\bar{u}_{\epsilon_j,\tau_j}(s)),\dot{\bar{u}}_{\epsilon_j,\tau_j}(s)\rangle\,ds\right).
\label{(++)}
\end{equation}
From Lemma~\ref{lemma2}, since by construction $\bigcup(a_j^i,b_j^i)\subseteq A_j^\eta\cap(s_j,t_j)$, we conclude
\begin{equation*}
F(t,u_-(t))-F(t,u_+(t))\geq\displaystyle\sum_{i=1}^m\mathop{\lim\sup}_{j\to+\infty}\int_{a_j^i}^{b_j^i}\|\nabla_x F(s,\bar{u}_{\epsilon_j,\tau_j}(s))\|\|\dot{\bar{u}}_{\epsilon_j,\tau_j}(s)\|\,ds.
\end{equation*}
Combining (\ref{(++)}) with Theorem~\ref{thm2.4} {(v)} we infer that
\begin{equation}
\begin{split}
F(t,u_-(t))-F(t,u_+(t))&\geq\displaystyle\sum_{i=1}^m\mathop{\lim\sup}_{j\to+\infty}\int_{a_j^i}^{b_j^i}\|\nabla_x F_{\tau_j}(s,\bar{u}_{\epsilon_j,\tau_j}(s))\|\|\dot{\bar{u}}_{\epsilon_j,\tau_j}(s)\|\,ds\\
&\geq\displaystyle\sum_{i=1}^m c_t(u_\eta^i, v_\eta^i),
\end{split}
\label{(+++)}
\end{equation}
where
\begin{equation*}
u_\eta^i=\displaystyle\lim_{j\to+\infty}\bar{u}_{\epsilon_j,\tau_j}(a_j^i) \text{ and }v_\eta^i=\displaystyle\lim_{j\to+\infty}\bar{u}_{\epsilon_j,\tau_j}(b_j^i),
\end{equation*}
whose existence is intended to be up to a subsequence of $j$. We note that these functions as well as the construction contained in this proof depend on $\eta$.

We now may fix a sequence of $\eta \to 0$ such that all the $u^i_{\eta}$'s and $v^i_\eta$'s converge, and the limit
\[
\lim_{\eta\to 0^+}\sum_{i=1}^m c_t(u_\eta^i, v_\eta^i)
\]
exists. Since by construction ${\rm dist}(u_\eta^i,\mathcal{C}(t))\leq\eta$, ${\rm dist}(v_\eta^i,\mathcal{C}(t))\leq\eta$ and $\|u_\eta^{i+1}-v_\eta^i\|\leq2\eta$, it follows that $u_\eta^{i+1}$ and $v_\eta^i$ have the same limit $u^{i+1}\in \mathcal{C}(t)$ as $\eta\to0^+$. Taking the limit in \eqref{(+++)} and using Theorem~\ref{thm2.4}, (vi) and (iv), we get
\begin{equation*}
F(t,u_-(t))-F(t,u_+(t))\geq \displaystyle\sum_{i=0}^{m-1}c_t(u^i,u^{i+1})\geq c_t(u^0,u^m)=c_t(u_-(t),u_+(t)),
\end{equation*}
since by (1) and (2) in the Claim it must be $u^0=u_-(t)$ and $u^m=u_+(t)$. Finally, again by Theorem~\ref{thm2.4} (ii), $c_t(u_+(t),u_-(t))=c_t(u_-(t),u_+(t))$ and this concludes the proof.
\endproof

The following proposition shows that, actually, \eqref{lowbounddiss} is an equality.

\begin{prop}\label{costotrovato}
Under the assumptions of {\rm Proposition \ref{taupiccolo}}, it holds
\begin{equation}
F(t,u_-(t))-F(t,u_+(t))=c_t(u_+(t),u_-(t)),\text{ for every }\,t\in[0,T].
\end{equation}
\end{prop}

\proof
In view of Proposition~\ref{proposizione1}, it will suffice to show the converse inequality
\begin{equation*}
F(t,u_-(t))-F(t,u_+(t)) \leq c_t(u_+(t),u_-(t)),\quad \text{ for every }\, t\in J.
\end{equation*}
For this, we take an arbitrary $\theta\in{\rm AC}([0,1];X)$ such that $\theta(0)=u_+(t)$, $\theta(1)=u_-(t)$ and we observe that
\begin{equation*}
\begin{split}
F(t,u_-(t))-F(t,u_+(t))&=F(t, \theta(1))-F(t,\theta(0))\\
&= \int_0^1 \langle \nabla_x F(t,\theta(s)), \dot{\theta}(s)\rangle\,ds\\
&\leq \int_0^1  \|\nabla_xF(t,\theta(s))\|\|\dot{\theta}(s)\|\,ds,
\end{split}
\end{equation*}
from which, taking the infimum on the right-hand side over the class of the admissible curves $\theta$, the thesis follows.
\endproof

We conclude with the following theorem that summarizes the results of this section and characterizes $u$ as a Balanced Viscosity solution of the problem
\begin{equation*}
\nabla_x F(t,u(t))=0\,\text{ in $X$}\quad\text{ for a.e. $t\in[0,T]$.}
\end{equation*}

\begin{thm}\label{main1}
Assume that {\bf (F0)-(F5)} hold and let $\tilde{u}_{\epsilon,\tau}: [0,T]\to X$ be the piecewise constant functions defined in \eqref{piecewiseconstant}, interpolating the discrete solutions of the implicit-time minimization scheme \eqref{scheme}, with $\tilde{u}_{\epsilon,\tau}(0)=u^0$. Let $c_t$ be the cost function defined in {\rm(\ref{costfun})}. Then all the sequences $(\varepsilon_j,\tau_j)_{j\in\mathbb{N}}$ satisfying $(\varepsilon_j,\tau_j)\to0$ and $\varepsilon_j/\tau_j\to+\infty$, as $j\to+\infty$, admit a subsequence {\rm(}still denoted by $(\varepsilon_j,\tau_j)${\rm)} and a limit function $u\colon[0,T]\mapsto X$ such that $\tilde{u}_{\epsilon_j,\tau_j}(t)\to u(t)$ for all $t \in [0, T]$. Moreover, $u$ satisfies the following properties:
\begin{enumerate}
\item[\rm(i)] $u$ is regulated;
\item[\rm(ii)] it holds
\begin{equation*}
\nabla_x F(t,u_+(t))=\nabla_xF(t, u_-(t))=0 \quad \text{in $X$ for every $t\in [0,T]$};
\end{equation*}
\item[\rm(iii)] $u$ fulfills the energy balance 
\begin{equation}
\begin{split}
F(t,u_+(t)) &+ \sum_{r\in J\cap[s,t]} c_r(u_-(r),u_+(r))=F(s,u_-(s))\\
    &+\int_s^t \partial_r F(r, u(r))\,dr,\quad \text{ for every }\,0\leq s\leq t\leq T.
\end{split}
\label{ebalance}
\end{equation}
\end{enumerate}
\label{bvthm}
\end{thm}

\proof
The result follows by combining Theorems \ref{compactnessthm} and \ref{genbalance} with \eqref{mu} and Proposition \ref{costotrovato}.
\endproof

\subsection{The critical regime $\displaystyle\lim_{j\to+\infty}\frac{\epsilon_j}{\tau_j}=\lambda\in(0,+\infty)$.}\label{casofinito}
Throughout this subsection we assume that the ratio $\frac{\epsilon_j}{\tau_j}$ tends to a strictly positive, {\it finite} limit $\lambda$ and characterise the dissipation cost as a discrete crease energy. We will prove our results assuming, to spare some notation, directly the equality
\begin{equation}\label{semplifico}
\frac{\epsilon_j}{\tau_j}=\lambda
\end{equation}
for all $j$. This is indeed not restrictive, since $\sum_{k \in \mathcal{K}_\tau}\|u^k-u^{k-1}\|^2$ is uniformly bounded along \eqref{scheme}, so that the term
\[
\left(\frac{\epsilon}{2\tau}-\frac{\lambda}2\right)\sum_{k \in \mathcal{K}_\tau}\|u^k-u^{k-1}\|^2
\]
is in general only a uniformly small remainder that we prefer to neglect.

We introduce, following~\cite[Section~3.4]{MinSav}, a parametric residual stability function $\mathcal{R}_\lambda$, $\lambda>0$. The name is motivated by the fact that $\mathcal{R}_\lambda$ provides a measure of the failure of the \emph{stability condition} for $(t,u)\in[0,T]\times X$, namely
\begin{equation*}
F(t,u)\leq F(t,v) +\frac{\lambda}{2}\|u-v\|^2,\quad\text{ for every $v\in X$}
\end{equation*}
(see \eqref{residualstability} below).

It can be defined as the difference between the energy $F(t,u)$ and its Moreau-Yosida regularization, as follows.

\begin{defn}
For every $t\in[0,T]$, $u\in X$ and $\lambda>0$, the \emph{residual stability function} is defined by
\begin{equation}
\mathcal{R}_\lambda(t,u):=F(t,u)-\displaystyle\min_{v\in X}\left\{F(t,v)+\frac{\lambda}{2}\|v-u\|^2\right\}.
\label{residual}
\end{equation} 
\end{defn}
\noindent
Note that $\mathcal{R}_\lambda(t,u)\geq0$ for every $t\in[0,T]$ and $u\in X$. Moreover,
\begin{equation}
\mathcal{R}_\lambda(t,u)=0\, \text{ if and only if }\, F(t,u)\leq F(t,v) +\frac{\lambda}{2}\|u-v\|^2,\quad\text{ for every $v\in X$,}
\label{residualstability}
\end{equation}
as it can be immediately checked from $\mathcal{R}_\lambda(t,u)=\displaystyle\max_{v\in X}\left\{F(t,u)-F(t,v)-\frac{\lambda}{2}\|v-u\|^2\right\}$. 

With the following proposition we provide some topological properties of the residual function $\mathcal{R}_\lambda$. 
\begin{prop}
Let $\mathcal{R}_\lambda$ be defined as in \eqref{residual}. Then the following properties hold:
\begin{enumerate}
\item[\rm(i)] $\mathcal{R}_\lambda(\cdot, u)$ is Lipschitz continuous uniformly on the compact subsets of $X$; i.e., for any $K\subset X$ compact, there exists $L_K>0$ such that for all $s,t\in[0,T]$ and $u\in K$ it holds
\begin{equation*}
|\mathcal{R}_\lambda(t, u)-\mathcal{R}_\lambda(s, u)|\leq L_K |t-s|;
\end{equation*}
\item[\rm(ii)] $\mathcal{R}_\lambda(t, \cdot)$ is continuous on $X$, for every $t\in[0,T]$;
\item[\rm(iii)] If we set 
\begin{equation}\label{zeroset}
\mathcal{Z}_\lambda(t):=\Bigl\{u\in X:\, \mathcal{R}_\lambda(t,u)=0\Bigr\}\,,
\end{equation}
then $\mathcal{Z}_\lambda(t)\subseteq \mathcal{C}(t)$ for every $t\in[0,T]$, where $\mathcal{C}(t)$ is defined as in \eqref{criticalset}. In particular, $\#\mathcal{Z}_\lambda(t)$ is locally  finite.
\end{enumerate}
\label{contres}
\end{prop}

\proof
(i) We fix $K\subset X$ compact subset and assume that $u\in K$. Let $v_t\in K$ be such that $v_t\in\displaystyle\mathop{\rm argmin}_{v\in X}\left\{F(t,v)+\frac{\lambda}{2}\|v-u\|^2\right\}$. Then we have
\[
\mathcal{R}_\lambda(t, u)= F(t,u)-F(t,v_t)-\frac{\lambda}{2}\|v_t-u\|^2,
\]
which combined with $F(t,v_t)\leq F(t,u)$ and {\bf (F0)-(F1)} gives
\[
|F(t,v_t)-F(s,v_t)|\leq L|t-s|,
\]
where $L=L_K$ is independent of $t$. Now,
\[
\begin{split}
\mathcal{R}_\lambda(t,u)&=(F(t,u)-F(t,v_t)-F(s,u)+F(s,v_t))\\
&+F(s,u)-F(s,v_t)-\frac{\lambda}{2}\|v_t-u\|^2\\
&\geq -2L|s-t|+\mathcal{R}_\lambda(s,u).
\end{split}
\]
Thus, 
\[
\mathcal{R}_\lambda(s,u)-\mathcal{R}_\lambda(t,u)\leq 2L|t-s|,
\]
from which the thesis follows exchanging the roles of $s$ and $t$. 

(ii) The proof easily follows by a standard $\Gamma$\hbox{-}convergence argument (see, e.g.,~\cite{GCB}). Let us fix $u\in X$ and let $\{u_n\}_{n\in\mathbb{N}}\subset X$ be such that $\|u_n-u\|\to0$ as $n\to+\infty$. Since $F(t,\cdot)$ is continuous from assumption {\bf (F0)}, it will suffice to show that
\begin{equation}
\displaystyle\min_{v\in X}\left\{F(t,v)+\frac{\lambda}{2}\|v-u_n\|^2\right\}\to\displaystyle\min_{v\in X}\left\{F(t,v)+\frac{\lambda}{2}\|v-u\|^2\right\},
\label{minconvergence}
\end{equation}
as $n\to+\infty$. Now, setting $G_n(t,v):=F(t,v)+\frac{\lambda}{2}\|v-u_n\|^2$ and $G(t,v):=F(t,v)+\frac{\lambda}{2}\|v-u\|^2$, we have that $G_n\to G$ for any $(t,v)$ uniformly as $n\to+\infty$. This implies that $G_n$ $\Gamma$\hbox{-}converge to $G$, so that the convergence of the minima \eqref{minconvergence} follows.

(iii) Let $u\in \mathcal{Z}_\lambda(t)$. By definition, we have
\begin{equation*}
F(t,u)=\displaystyle\min_{v\in X}\left\{F(t,v)+\frac{\lambda}{2}\|v-u\|^2\right\},
\end{equation*}
so that
\begin{equation*}
\nabla_x\left[F(t,v)+\frac{\lambda}{2}\|v-u\|^2\right]_{|v=u}=0.
\end{equation*}
Thus, $\nabla_xF(t,u)=0$, and therefore $u\in \mathcal{C}(t)$.
\endproof

\begin{defn}
We define the \emph{transition cost} $c^\lambda$ as
\begin{equation}
c^\lambda(t,u,v):=\displaystyle\inf_{N\in\mathbb{N}}\left\{\sum_{i=0}^{N-1}\frac{\lambda}{2}\|w^i-w^{i+1}\|^2+\sum_{i=0}^{N}\mathcal{R}_\lambda(t,w^i):\, W:\mathbb{N}_0\to X,\, w^0=u,\, w^N=v\right\},
\label{cost}
\end{equation}
where $\mathbb{N}_0=\mathbb{N}\cup\{0\}$ and $W$ denotes the family $(w^i)_{i=0,\dots,N}$.
\end{defn}

\begin{oss}
The transition cost \eqref{cost} satisfies the following elementary properties:
\begin{enumerate}
\item[(i)] $c^\lambda(t,u,v)\leq c^\lambda(t,u,w)+c^\lambda(t,w,v)$, \quad $\forall u,v,w\in X$;
\item[(ii)] $c^\lambda(t,u,v)=c^\lambda(t,v,u)$, \quad $\forall u,v\in X$.
\end{enumerate}
Furthermore $c^\lambda(t,u,u)=0 \iff u \in \mathcal{Z}_\lambda(t)$.
\end{oss}

Now we prove a semicontinuity property which will be fundamental in the sequel.

\begin{thm}
Let $(m_j)_j$ be a sequence of positive integers, and $I_j=\{0,1,\dots,m_j\}$. For every $j$ consider $T_j=(t_j^i)_{i\in I_j}$ satisfying 
\[
t_j^i<t_j^{i+1}, \quad  t_j^0\to t, \quad t_j^{m_j}\to t,
\]
and $W_j=(w_j^i)_{i\in I_j}$ such that 
\[
W_j\subseteq B_M, \quad w_j^0\to u, \quad w_j^{m_j}\to v\,.
\]
Assume \eqref{semplifico} and, for $ \mathcal{Z}_\lambda(t)$ as in \eqref{zeroset}, suppose that $u,v\in \mathcal{Z}_\lambda(t)$,  with $u\neq v$. Then
\begin{equation}
\displaystyle\mathop{\lim\inf}_{j\to+\infty}\left[\sum_{i=0}^{m_j - 1}\frac{\lambda}{2}\|w_j^i-w_j^{i+1}\|^2+\mathcal{R}_\lambda(t_j^{i+1}, w_j^i)\right]\geq c^\lambda(t,u,v).
\label{stimata}
\end{equation}
\label{semicont}
\end{thm}

\proof
It is sufficient to treat the case $\displaystyle\lim_{j\to+\infty}m_j=+\infty$, the proof being similar (and simpler) if $\displaystyle\mathop{\lim\inf}_{j\to+\infty}m_j<+\infty$. It is not restrictive to assume that the liminf in the left-hand side of \eqref{stimata} is a limit, and that it is finite (otherwise the result is trivial).

Let $\eta>0$ be fixed. Correspondingly, we define the set
\begin{equation*}
G_j^\eta:=\Bigl\{i\in I_j:\, {\rm dist}(w_j^i,\mathcal{Z}_\lambda(t))\geq\eta\Bigr\}.
\end{equation*}
First, we prove that there exists a positive constant $C_\eta$ such that 
\begin{equation}\label{boundG}
\#G_j^\eta\leq C_\eta<+\infty,\quad \text{ for every $j$.} 
\end{equation}
For this, by the continuity of function $\mathcal{R}_\lambda(t,\cdot)$ (Proposition~\ref{contres} (ii)), there exists a constant $\gamma_\eta$ such that $\mathcal{R}_\lambda(t,v)\geq\gamma_\eta$ for any $v \in B_M$ satisfying ${\rm dist}(v,\mathcal{Z}_\lambda(t))\geq\eta$. Thus, for any $i\in G_j^\eta$, we have $\mathcal{R}_\lambda(t,w_j^i)\geq\gamma_\eta$. On the other hand, we may use the Lipschitz continuity of $\mathcal{R}_\lambda(\cdot,u)$ (Proposition~\ref{contres} (i)) to find
\begin{equation*}
|\mathcal{R}_\lambda(t_j^{i+1},w_j^i)-\mathcal{R}_\lambda(t,w_j^i)|\leq L|t_j^{i+1}-t|\leq L \max\{|t_j^0-t|,|t_j^{m_j}-t|\}=o(1),
\end{equation*}
as $j\to+\infty$. We then have
\begin{equation*}
\mathcal{R}_\lambda(t_j^{i+1},w_j^i)\geq \frac{1}{2}\gamma_\eta,
\end{equation*}
for $j$ large enough. Since $\displaystyle\sum_{i=0}^{m_j - 1}\mathcal{R}_\lambda(t_j^{i+1}, w_j^i)\leq C$ holds,  as the limit in \eqref{stimata} is finite, the above inequality forces $\#G_j^\eta$ to be bounded from above by a constant depending only on $\eta$. In particular, $\#G_j^\eta$ is finite, as desired.

Now, in order to prove \eqref{stimata}, we argue by induction on the index $i$. We define
\begin{equation*}
i_j^{1,-}:=\max\Bigl\{i\in I_j:\, \|w_j^i-u\|\leq\eta\Bigr\}\,.
\end{equation*}
Since $ \mathcal{Z}_\lambda(t)$ consists of isolated points by {\bf(F3)} and Proposition~\ref{contres}(iii), and $W_j$ is confined to the bounded set $B_M$, by the condition $w_j^{m_j}\to v$ it holds $i_j^{1,-}<m_j$, provided that $\eta$ is suitably chosen. We then set
\begin{equation}
i_j^{1,+}:=\min\left\{i>i_j^{1,-}:\, {\rm dist}(w_j^i,\mathcal{Z}_\lambda(t))\leq\eta \right\}.
\label{i+}
\end{equation}
Note that the minimum in (\ref{i+}) is well-defined, since the set contains at least $m_j$ and is therefore nonempty. Moreover, if $z\in \mathcal{Z}_\lambda(t)$ is such that $\left\|w_j^{i_j^{1,+}}-z\right\|= {\rm dist}\left(w_j^{i_j^{1,+}},\mathcal{Z}_\lambda(t)\right)$, then necessarily $z\neq u$. Now, if $z=v$ then the proof stops here, otherwise the inductive argument goes on.

Let $z_{k-1}\in \mathcal{Z}_\lambda(t)$ be such that $\left\|w_j^{i^{k-1,+}}-z_{k-1}\right\|= {\rm dist}\left(w_j^{i^{k-1,+}},\mathcal{Z}_\lambda(t)\right)$, for some $k\geq2$; by induction it results $z_{k-1}\neq v$. We set
\begin{equation*}
i_j^{k,-}:=\max\left\{i\geq i_j^{k-1,+}:\, \|w_j^i-z_{k-1}\|\leq\eta\right\},
\end{equation*}
and
\begin{equation*}
i_j^{k,+}:=\min\left\{i> i_j^{k,-}:\, {\rm dist}(w_j^{i},\mathcal{Z}_\lambda(t))\leq\eta\right\},
\end{equation*}
where, as before, the minimum is well-defined since $m_j$ belongs to the set. The argument stops when the point in $\mathcal{Z}_\lambda(t)$ with minimal distance from $w_j^{i_j^{k,+}}$ coincides with $v$, and this happens after at most  $\#\left(\mathcal{Z}_\lambda(t)\cap B_M\right)=:N_t$ steps, with $N_t<+\infty$ by {\bf(F3)} and Proposition~\ref{contres}(iii).

We may assume, without loss of generality, that the required number of steps in order to conclude the inductive procedure is exactly $N_t$.
Now we note that 
\begin{equation}
M_j^t:=\displaystyle\sum_{k=1}^{N_t}(i_j^{k,+}-i_j^{k,-})\leq C_\eta+N_t<+\infty.
\label{C*}
\end{equation}
Indeed, by construction, each $i$ such that $i_j^{k,-}<i<i_j^{k,+}$ satisfies ${\rm dist}(w_j^i, \mathcal{Z}_\lambda(t))>\eta$ and then $i\in G_j^\eta$. With this, \eqref{C*} follows from \eqref{boundG} Moreover, it holds by construction
\begin{equation}
\left\|w_j^{i_j^{k,+}}-w_j^{i_j^{k+1,-}}\right\|\leq 2\eta,\,\text{ for every }\, k.
\label{C**}
\end{equation}

Setting $I_j^t:=\bigcup_{k=1}^{N_t}\{i_j^{k,-},\dots, i_j^{k,+}\}$, we consider the set of positive integers $L_j^t:=\{0\}\cup I_j^t \cup\{M_j^t+1\}$. We put $w_j^0=u$, $w_j^{M_j^t+1}=v$.
For $l\in\{1,\dots, M_j^t\}$ we set $w_j^l=w_j^i$, where $i$ is the $l$-th point in $I_j^t$; analogously, $t_j^l=t_j^i$. We now claim that 
\begin{equation}\label{stimette}
\begin{split}
\displaystyle\sum_{i=1}^{m_j-1}\left[\frac{\lambda}{2}\|w_j^i-w_j^{i+1}\|^2+\mathcal{R}_\lambda(t_j^{i+1},w_j^i)\right]&\geq\displaystyle\sum_{l=0}^{M_j^t}\left[\frac{\lambda}{2}\|w_j^l-w_j^{l+1}\|^2+\mathcal{R}_\lambda(t_j^{l+1},w_j^l)\right]\\
&-2\lambda\eta^2(N_t+1)\,.
\end{split}
\end{equation}
where we used the convention that $t_j^{M_j^t+1}=t_j^{M_j^t}+\frac{1}{j}$.  Indeed, by construction we have that $\|w_j^1-u\|\leq\eta$ and $\|w_j^{M_j^t}-v\|\leq\eta$. With this and \eqref{C**}, \eqref{stimette} follows by noting that the only terms of the sum in the right-hand side which do not appear in the left-hand side are those of the form
\[
\frac{\lambda}{2}\left\|u-w_j^{i_j^{1,-}}\right\|^2,\quad \frac{\lambda}{2}\left\|w_j^{i_j^{k,+}}-w_j^{i_j^{k+1,-}}\right\|^2, \quad \frac{\lambda}{2}\left\|w_j^{M_j^t}-v\right\|^2 .
\]
If we further add and subtract in the right-hand side of \eqref{stimette} the term $\mathcal{R}_\lambda(t,w_j^l)$ and noting that $\mathcal{R}_\lambda(t,w_j^{M_j^t+1})=\mathcal{R}_\lambda(t,v)=0$, we obtain
\begin{equation*}
\begin{split}
\displaystyle\sum_{i=1}^{m_j-1}\left[\frac{\lambda}{2}\|w_j^i-w_j^{i+1}\|^2+\mathcal{R}_\lambda(t_j^{i+1},w_j^i)\right]&\geq\displaystyle\sum_{l=0}^{M_j^t}\frac{\lambda}{2}\|w_j^l-w_j^{l+1}\|^2+\displaystyle\sum_{l=0}^{M_j^t+1}\mathcal{R}_\lambda(t, w_j^l)\\
&-2\lambda\eta^2(N_t+1)-\zeta_j,
\end{split}
\end{equation*}
where 
\begin{equation*}
\zeta_j=\left|\sum_{l=0}^{M_j^t}\left[\mathcal{R}_\lambda(t_j^{l+1}, w_j^l)-\mathcal{R}_\lambda(t,w_j^l)\right]\right|=o(1),\quad\text{ as }j\to+\infty,
\end{equation*}
from the Lipschitz continuity of $\mathcal{R}_\lambda(\cdot,w)$ and the equi-boundedness of $M_j^t$ (with a constant depending on $\eta$) ensured by \eqref{C*}.

Now, since $(w_j^l)_{l=0,\dots, M_j^t+1}$ is an admissible test function in the minimization problem \eqref{cost} defining $c^\lambda(t,u,v)$, we conclude that
\begin{equation*}
\displaystyle\mathop{\lim\inf}_{j\to+\infty}\left[\sum_{i=0}^{m_j - 1}\frac{\lambda}{2}\|w_j^i-w_j^{i+1}\|^2+\mathcal{R}_\lambda(t_j^{i+1}, w_j^i)\right]\geq c^\lambda(t,u,v)-2\lambda\eta^2(N_t+1).
\end{equation*}
Finally, since $N_t=\#(\mathcal{Z}_\lambda(t)\cap B_M)$ is independent of $\eta$, letting $\eta\to0$ in the latter inequality we get the thesis.
\endproof

We can now state and prove a complete characterization in terms of stability  and energy balance of the limit behavior of \eqref{scheme}, whenever the ratio $\frac{\epsilon}{\tau}$ tends to a finite limit $\lambda$.

\begin{thm}\label{main2}
Assume that {\bf (F0)-(F5)} hold and let $\tilde{u}_{\epsilon,\tau}: [0,T]\to X$ be the piecewise constant functions defined in \eqref{piecewiseconstant}, interpolating the discrete solutions of the implicit-time minimization scheme \eqref{scheme}, with $\tilde{u}_{\epsilon,\tau}(0)=u^0$. Let $\lambda>0$ be a real number and $c^\lambda$ be the transition cost function defined in {\rm(\ref{cost})}. Then all the sequences $(\varepsilon_j,\tau_j)_{j\in\mathbb{N}}$ satisfying $(\varepsilon_j,\tau_j)\to0$ and $\varepsilon_j/\tau_j\to\lambda$, as $j\to+\infty$, admit a subsequence {\rm(}still denoted by $(\varepsilon_j,\tau_j)${\rm)} and a limit function $u\colon[0,T]\mapsto X$ such that $\tilde{u}_{\epsilon_j,\tau_j}(t)\to u(t)$ for all $t \in [0, T]$. Moreover, $u$ satisfies the following properties:
\begin{enumerate}
\item[\rm(i)] $u$ is regulated;
\item[\rm(ii)] \emph{(Stability)} $u$ satisfies the stability conditions
\begin{equation*}
F(t, u_+(t))\le F(t, v)+\frac{\lambda}2\|v-u_+(t)\|^2
\end{equation*}
and
\[
F(t, u_-(t))\le F(t, v)+\frac{\lambda}2\|v-u_-(t)\|^2
\]
for all $t\in [0,T]$ and all $v\in X$.
\item[\rm(iii)] \emph{(Energy balance)} $u$ fulfills the energy balance
\begin{equation}\label{ebalance2}
\begin{split}
F(t,u_+(t))&+\displaystyle\sum_{r\in J\cap[s,t]} c^\lambda(r, u_+(r),u_-(r))\\
&= F(s,u_-(s))+\int_s^t\partial_r F(r,u(r))\,dr,\quad \text{ for every }\,0\leq s\leq t\leq T.
\end{split}
\end{equation}
\end{enumerate}
\end{thm}

\proof
In view of Theorem \ref{compactnessthm}, Proposition \ref{improve}, Theorem \ref{genbalance},  and \eqref{mu}, we are left to show that
\begin{equation}\label{ziel}
F(t,u_-(t))-F(t,u_+(t))=c^\lambda(t, u_-(t),u_+(t)),\quad \text{ for every }\, t\in J\,,
\end{equation}
where $J$ denotes as usual the jump set of $u$. Let $\eta>0$ and $t\in J$ be fixed. Then by the definition of $c^\lambda$, there exists a positive integer $N_\eta$ and a family of functions $U=(u^i)_{i=0,\dots, N_\eta}$, with $u^0=u_-(t)$ and $u^{N_\eta}=u_+(t)$, such that
\begin{equation}
\sum_{i=0}^{N_\eta-1}\frac{\lambda}{2}\|u^i-u^{i+1}\|^2+\mathcal{R}_\lambda(t,u^i)\leq c^\lambda(t,u_-(t),u_+(t))-\eta.
\label{equality}
\end{equation}
Notice that  the term $\mathcal{R}_\lambda(t,u^{N_\eta})$ does not contribute to the sum above, since it is $0$ because of the stability condition (ii).
Now, from the definition of $\mathcal{R}_\lambda$ we have
\[
\begin{split}
\frac{\lambda}{2}\|u^i-u^{i+1}\|^2+\mathcal{R}_\lambda(t,u^i)&\geq \frac{\lambda}{2}\|u^i-u^{i+1}\|^2+F(t,u^i)-F(t,u^{i+1})-\frac{\lambda}{2}\|u^i-u^{i+1}\|^2\\
&=F(t,u^i)-F(t,u^{i+1})
\end{split}
\]
for all $i=0,1,\dots,N_\eta -1$, so that
\begin{equation*}
\sum_{i=0}^{N_\eta-1}\frac{\lambda}{2}\|u^i-u^{i+1}\|^2+\mathcal{R}_\lambda(t,u^i)\geq \sum_{i=0}^{N_\eta-1} \Bigl[F(t,u^i)-F(t,u^{i+1})\Bigr]\\
= F(t,u_-(t))-F(t,u_+(t)).
\end{equation*}
Combining this estimate with \eqref{equality} and letting $\eta\to0$ we recover the first inequality
\begin{equation}\label{ziel12}
F(t,u_-(t))-F(t,u_+(t))\leq c^\lambda(t, u_-(t),u_+(t)).
\end{equation}

To prove the converse, by a {diagonal} argument and \eqref{limite}, we may find two sequences $t_j\searrow t$ and $s_j\nearrow t$ such that
\begin{equation*}
\mathop{\lim\inf}_{j\to+\infty} \Bigl[\|\tilde{u}_{\epsilon_j,\tau_j}(s_j)-u_-(t)\|+\|\tilde{u}_{\epsilon_j,\tau_j}(t_j)-u_+(t)\|\Bigr]=0,
\end{equation*}
so that
\begin{equation*}
\mathop{\lim\inf}_{j\to+\infty} \Bigl[|F(s_j,\tilde{u}_{\epsilon_j,\tau_j}(s_j))-F(t,u_-(t))|+|F(t_j,\tilde{u}_{\epsilon_j,\tau_j}(t_j))-F(t,u_+(t))|\Bigr]=0,
\end{equation*}
and then 
\begin{equation}
F(t,u_-(t))-F(t,u_+(t))\geq \mathop{\lim\inf}_{j\to+\infty}\left(F(s_j,\tilde{u}_{\epsilon_j,\tau_j}(s_j))-F(t_j,\tilde{u}_{\epsilon_j,\tau_j}(t_j))\right).
\label{star}
\end{equation}
As already remarked in other proofs, it is not restrictive to assume that $s_j$ and $t_j$ are nodes in $\Pi_j$. Setting $m_j+1:=\#\{t_j^i\in \Pi_j:\, s_j\leq t_j^i\leq t_j\}$, $t_j^0=s_j$ and $t_j^{m_j}=t_j$ we have
\begin{equation*}
\begin{split}
F(s_j,\tilde{u}_{\epsilon_j,\tau_j}(s_j))-F(t_j,\tilde{u}_{\epsilon_j,\tau_j}(t_j))&=\sum_{i=0}^{m_j-1}F(t_j^i,\tilde{u}_{\epsilon_j,\tau_j}(t_j^i))-F(t_j^{i+1},\tilde{u}_{\epsilon_j,\tau_j}(t_j^{i+1}))\\
&\geq -L(t_j-s_j)
\\&+\sum_{i=0}^{m_j-1}\left[F(t_j^{i+1},\tilde{u}_{\epsilon_j,\tau_j}(t_j^i))-F(t_j^{i+1},\tilde{u}_{\epsilon_j,\tau_j}(t_j^{i+1}))\right],
\end{split}
\end{equation*}
where $L$ is the Lipschitz constant of $F(\cdot,u)$ uniform for $u\in B_M$. From the definition of $\mathcal{R}_\lambda$, \eqref{semplifico} and the minimality of $\tilde{u}_{\epsilon_j,\tau_j}(t^{i+1}_j)$ it follows that
\begin{equation*}
\begin{split}
&F(s_j,\tilde{u}_{\epsilon_j,\tau_j}(s_j))-F(t_j,\tilde{u}_{\epsilon_j,\tau_j}(t_j)) + L(t_j-s_j)\\
&\geq\sum_{i=0}^{m_j-1}\left[\frac{\lambda}{2}\|\tilde{u}_j(t_j^i)-\tilde{u}_{\epsilon_j,\tau_j}(t_j^{i+1})\|^2
+\mathcal{R}_\lambda(t_j^{i+1},\tilde{u}_{\epsilon_j,\tau_j}(t_j^{i+1}))\right]. 
\end{split}
\end{equation*}
Since $(t_j-s_j)\to0$, as a consequence of (\ref{star}) and Theorem~\ref{semicont} we deduce the converse inequality
\begin{equation*}
F(t,u_-(t))-F(t,u_+(t))\geq c^\lambda(t, u_-(t),u_+(t))\,.
\end{equation*}  
Combining with \eqref{ziel12} we get \eqref{ziel} and the proof is complete.
\endproof

\end{document}